\documentclass{article}
\usepackage{graphicx} % Required for inserting images
\usepackage[affil-it]{authblk}
\usepackage[utf8]{inputenc}
\usepackage[english]{babel}
\usepackage{titlesec}
\usepackage[margin=0.5in]{geometry}
\usepackage{mathtools}
\usepackage{enumerate}
\usepackage{enumitem}
\usepackage{lineno} % Add line numbers
\usepackage{adjustbox}
\usepackage{multirow}
\usepackage{chngcntr}
\usepackage{titlesec}
\usepackage[utf8]{inputenc}
\usepackage{mathrsfs}
\usepackage{lscape}
\usepackage{cite}
\usepackage{cite}

\usepackage{enumitem}
\usepackage{subfigure}
\usepackage{subcaption}

\usepackage[hyphens]{url}
\usepackage{hyperref}
\usepackage{float}
\usepackage[section]{placeins}
\usepackage[dvipsnames]{xcolor}
\usepackage{soul}
\usepackage{amsmath}
\usepackage[labelfont=bf,font=small]{caption}
\usepackage{bm}
\usepackage[title]{appendix}
\usepackage{graphicx}
\usepackage{tabularx}
\usepackage{float}
\usepackage{mathtools}
\usepackage{amsmath,amssymb,amsthm}
\usepackage{oubraces}
\usepackage{chngcntr}
\usepackage{mathrsfs}
\usepackage{multirow}
\usepackage{chngcntr}
\usepackage{parskip}
\usepackage{accents}
\usepackage{soul}
\usepackage{orcidlink}
\usepackage[normalem]{ulem}

 \usepackage{changepage}

\usepackage[dvipsnames]{xcolor}   %latex colors

\definecolor{dgn}{rgb}{0.0, 0.5, 0.0}

\newcommand{\allblack}{\color{black}{}}

\usepackage{lipsum}% http://ctan.org/pkg/lipsum
\makeatletter
\g@addto@macro{\endabstract}{\@setabstract}
\newcommand{\authorfootnotes}{\renewcommand\thefootnote{\@fnsymbol\c@footnote}}%
\makeatother
\begin{document}
\title{Cell Migration Boundary Motion in \emph{Drosophila} Egg Chambers: A combined Phase Field and Chemoattractant Model}
\author{Naghmeh Akhavan \orcidlink{0000-0002-9474-4486}$^{a}$, Alexander S. George$^{b}$, Ahmet Kaan Aydin\orcidlink{0009-0006-6278-2046}$^{d}$,  Michelle Starz-Gaiano\orcidlink{0000-0002-0855-0715}$^{c}$,  and Bradford E. Peercy \orcidlink{0000-0002-8597-2508}$^{d}\footnote{Corresponding author: Bradford E. Peercy {\it (bpeercy@umbc.edu)}}$\\
% \begin{flushleft}
{\it {\small $^{a}$  Department of Mathematics, University of Michigan, Ann Arbor, MI, USA}}

 {\it {\small $^{b}$  Dartmouth Cancer Center, One Medical Center Drive, Lebanon, NH, USA}}

 {\it {\small $^{c}$  Department of Biological Sciences, University of Maryland Baltimore County, MD, USA}}

{\it {\small $^{d}$  Department of Mathematics and Statistics, University of Maryland Baltimore County, MD, USA}}

% \end{flushleft}
%    \today
\vspace{-1.5cm}
    }
  \date{}
\maketitle
\vspace{1cm}

\section*{Abstract}

In the \emph{Drosophila melanogaster} egg chamber, collective migration of border cells toward the oocyte is guided by spatial chemoattractant gradients, yet the distribution of these signals within the tissue remains difficult to quantify experimentally. We develop a spatially resolved mathematical framework to examine how tissue geometry shapes local chemoattractant concentrations during migration. A phase-field model represents the egg chamber architecture and is coupled to a diffusion--reaction equation with spatially heterogeneous diffusivity, restricting free chemoattractant transport primarily to extracellular space while reducing diffusion within cellular domains. Chemoattractant is secreted from the oocyte boundary and degraded throughout the domain, producing geometry-dependent concentration patterns. We further couple this chemical field to a mechanical model of cluster migration through the Tangential Interface Migration (TIM) force, which links local chemoattractant cues with contact-mediated motion along neighboring nurse-cell interfaces. Simulations show that tissue architecture and signal localization jointly influence migration speed and directional persistence. In particular, geometric bottlenecks and intersections can flatten local gradients and slow migration. This framework provides a mechanistic approach for investigating how extracellular geometry, chemical signaling, and cell--cell interactions together regulate collective cell migration \emph{in vivo}.
\section{Introduction}

Border cell migration in the fruit fly \emph{Drosophila melanogaster} egg chamber offers a powerful system for studying collective cell migration in a confined, multicellular environment~\cite{rorth2009collective, montell2012group}.
During mid-oogenesis, a small cluster of epithelial cells, the border cells, detaches from the anterior follicular epithelium and migrates posteriorly through a crowded field of nurse cells to reach the oocyte~\cite{rorth2009collective, montell2012group, peercy2020clustered}. 
{\allblack Successful arrival at the oocyte has a direct developmental function: the border cells contribute to the formation and shaping of the micropyle, a specialized structure of the eggshell that permits sperm entry during fertilization \cite{montell2012group, horne2020drosophila}. 
This migration provides an experimentally accessible {\it in vivo} system in which collective cell cohesion, chemical guidance, and mechanical interactions with a confining tissue environment must be coordinated during development \cite{montell2012group, george2025traveling}.}
This movement is largely guided by spatial gradients of several different chemoattractants, including fruit fly homologs of Epidermal Growth Factor (EGF) and Platelet Derived Growth Factor and Vascular Endothelial Growth Factor (PDGF/VEGF)-family member ligands secreted by the oocyte~\cite{duchek2001guidance, mcdonald2003pvf1, duchek2001guidance2, dai2020tissue, montell2012group} and sensed by receptor tyrosine kinases on the surface of border cells~\cite{duchek2001guidance2, mcdonald2006multiple, duchek2001guidance}.  The  PDGF/VEGF ligand found in egg chambers is called PVF1~\cite{duchek2001guidance}.  Receptors at the border cell surface are tightly regulated to allow cells to respond to a dynamic range of ligand concentration and, upon activation, promote cytoskeletal changes and forward movement~\cite{wan2013guidance, boutet2023arfgap1, assaker2010spatial, zhou2022two, poukkula2011cell, bianco2007two, prasad2007cellular, george2025traveling}. {\allblack Notably, the molecular mechanisms that control border cell migration are well-conserved with those governing other types of developmental cell migrations \cite{aman2010cell, cheung2025collective} 
% (need citations doi: 10.1016/j.ydbio.2009.11.014 \&  DOI: 10.1038/s41580-025-00858-9).
Thus, detailed characterization of this process will likely reveal additional shared principals.}
% {\allred (cite Montell/Rorth reviews and Emery papers i will find)}

Although chemotactic signaling plays a central role in directing border cell migration~\cite{prasad2007cellular}, the actual distribution of chemoattractant is thought to be spatially heterogeneous and shaped by the geometry of the extracellular space~\cite{george2025chemotaxis}. {\allblack Because physical structure itself can also impact cell movements  ~\cite{dai2020tissue, aranjuez2016dynamic} , it is challenging to separately characterize the influence of physical and chemical guidance inputs.} Experimental evidence suggests that runways (bottlenecks) and intersections between nurse cells can locally flatten the chemoattractant gradient, often leading to slower or less directed migration~\cite{george2025chemotaxis}. These observations suggest that tissue architecture and spatial constraints significantly modulate signal availability and, in turn, the effectiveness of chemotactic cues.

% Directly measuring chemoattractant distributions \emph{in vivo} remains technically challenging due {to low signal levels, limited detection ability \textit{in vivo}, changes over time,} and limited spatial resolution~\cite{prasad2007cellular, sarris2015navigating}.
% % {(\allcblue maybe more references here!)}.
% While tagged ligands~\cite{dai2020tissue, george2025chemotaxis}
% {
% % (need to cite George and Dai here (for pvf1 and keren), respectively)
% offer some insight, their interpretation is complicated by photobleaching, potential kinetic changes due to tagging or overexpression, and the ability to monitor signals accurately over time and space}~\cite{dai2020tissue, george2025chemotaxis}. 
% % {(\allcblue maybe more references here! i think we are okay to just add Dai)}. 
% These limitations necessitate mathematical modeling to infer spatial gradient structures and evaluate how tissue geometry shapes signaling.

{\allblack Recent advances in fluorescence imaging, image analysis, and computational reconstruction have substantially improved the ability to quantify molecular distributions at subcellular resolution in developing tissues, including measurements of adhesion and cytoskeletal components in living embryos \cite{li2025ctransformer, yamamoto2026spatiotemporal, junyent2024first, soto2025spindle}. 
Nevertheless, quantitative measurement of low-level, diffusible extracellular guidance cues over space and time remains challenging, particularly in geometrically complex multicellular tissues such as the Drosophila egg chamber. Tagged ligands provide valuable information about chemoattractant localization \cite{george2025chemotaxis, dai2020tissue}, but quantitative interpretation can still be affected by signal intensity, photobleaching, temporal resolution, and potential changes in ligand kinetics associated with tagging or overexpression. Mathematical modeling therefore provides a complementary approach for examining spatial gradient structures that are difficult to resolve directly and for testing how tissue geometry may influence chemoattractant distribution.}

To address this, mathematical models have been developed to reconstruct chemoattractant dynamics and assess how spatial structure influences signal propagation~\cite{mekus2018effects, camley2016emergent, yue2018minimal, cai2016modeling, peercy2020clustered}. Earlier efforts, including compartmental and agent-based models~\cite{montell2012group, stonko2015mathematical}, often idealize tissue as a homogeneous medium or overlook signal exclusion from intracellular regions. More recent approaches leverage continuum phase-field methods to explicitly model extracellular space and cellular geometry, providing a natural framework for simulating chemically guided migration in biologically realistic environments. {\allblack More broadly, phase-field approaches have been increasingly applied to multicellular morphogenesis in other confined developmental systems, particularly the {\it Caenorhabditis elegans} embryo \cite{seirin2021extra, kuang2022computable, kuang2023morphosim}. 
Such approaches are particularly useful when cell shape, deformable interfaces, cell–cell contacts, and confinement need to be represented explicitly, and comparisons with coarse-grained descriptions have highlighted the additional cell-resolved information available from phase-field formulations \cite{guan2023comparison}. 
More recent extensions have sought to connect phase-field descriptions more directly with imaging data and to integrate multicellular mechanics with intercellular signaling \cite{sukekawa2025imaging, guan2026efficient}.}

In this study, we extend our previous phase-field model of border cell migration~\cite{akhavan2026phase} by developing a spatially resolved framework for chemoattractant distribution that captures key biophysical constraints observed \emph{in vivo}. Specifically, we incorporate:
(1) Boundary-localized secretion of chemoattractant from the oocyte, modeled using the gradient of the oocyte’s phase field variable;
(2) Spatially heterogeneous diffusion constrained by the geometry of cellular domains, allowing free diffusion only through extracellular corridors;
(3) A chemotactic force formulation that biases border cell movement in response to local gradients of chemoattractant concentration.

Through numerical simulations, we show that the spatial organization of the extracellular space strongly influences chemoattractant gradients, receptor activation, and migration behavior. Narrow intercellular gaps lead to anisotropic and spatially diminished
% {\allcpink is there some other word beside filtered we can use? i don't think it is very clear} {\allorange maybe diminished?}
gradients that can limit directional cues, especially at tissue intersections where multiple large cells meet. These results suggest that effective chemical guidance arises not solely from localized secretion, but also from the heterogeneous diffusion landscape shaped by tissue structure.

Importantly, we integrate this chemical framework with a mechanical model of migration that includes the Tangential Interface Migration (TIM) force~\cite{akhavan2026phase}, which captures contact-mediated propulsion along interfaces between the cluster and adjacent nurse cells. This allows us to examine how chemical gradients and mechanical forces interact—cooperatively or compensatorily—to govern the directional migration of the border cell cluster in a complex, multicellular environment.

% \section*{Materials and methods}

% Place figure captions after the first paragraph in which they are cited.
% \begin{figure}[!h]
% \caption{{\bf Bold the figure title.}
% Figure caption text here, please use this space for the figure panel descriptions instead of using subfigure commands. A: Lorem ipsum dolor sit amet. B: Consectetur adipiscing elit.}
% \label{fig1}
% \end{figure}

% Results and Discussion can be combined.
\section{Results}
We apply our extended phase-field model to examine how spatial chemoattractant dynamics and contact-mediated forces influence the collective migration of the border cell cluster in the \emph{Drosophila} egg chamber.

\subsection{Chemoattractant Diffusion and Concentration Across the Extracellular Space}

% We first simulate the spatially heterogeneous diffusion of the chemoattractant secreted by the oocyte into the extracellular domain of the egg chamber. Diffusion is restricted within cellular regions and primarily occurs through narrow gaps between cells (extracellular regions), mimicking the constrained extracellular space observed \emph{in vivo} \cite{george2025chemotaxis}. This behavior is encoded through a spatially varying diffusion coefficient \( D(\phi) \), which is dynamically modulated by the phase field variables \( \phi_i \) representing individual cells (Eq.~\eqref{eq pfm diff}). Specifically, \( D(\phi) \) is significantly reduced in regions where \( \phi_i \approx 1 \), corresponding to intracellular domains, and remains nonzero only in extracellular regions.

{\allblack We first examine how the cellular geometry of the egg chamber shapes the chemoattractant distribution prior to border-cell migration. The oocyte acts as the chemoattractant source, while the nurse cells, border-cell cluster, and oocyte contribute to the spatially varying diffusion landscape through their phase-field representations. According to Eq.~\eqref{eq pfm diff}, the effective diffusion coefficient is high in extracellular regions and strongly reduced within cellular domains. Consequently, chemoattractant transport occurs preferentially through the extracellular spaces separating neighboring cells.}
% and Figure~\ref{fig:AllenCahn}).

% \begin{figure}[H]
%     \centering
%     \includegraphics[width=0.65\linewidth]{Figures/1Ddiffsiont=16.pdf}
%     \caption{Evolution of the chemoattractant concentration $c(x,t)$, the  phase field variable $\phi(x,t)$, and the diffusion coefficient $D(\phi(x,t))$ over time.}
%     \label{fig:AllenCahn}
% \end{figure}

\begin{figure}[htp!]
    \centering
    \includegraphics[width=1\linewidth]{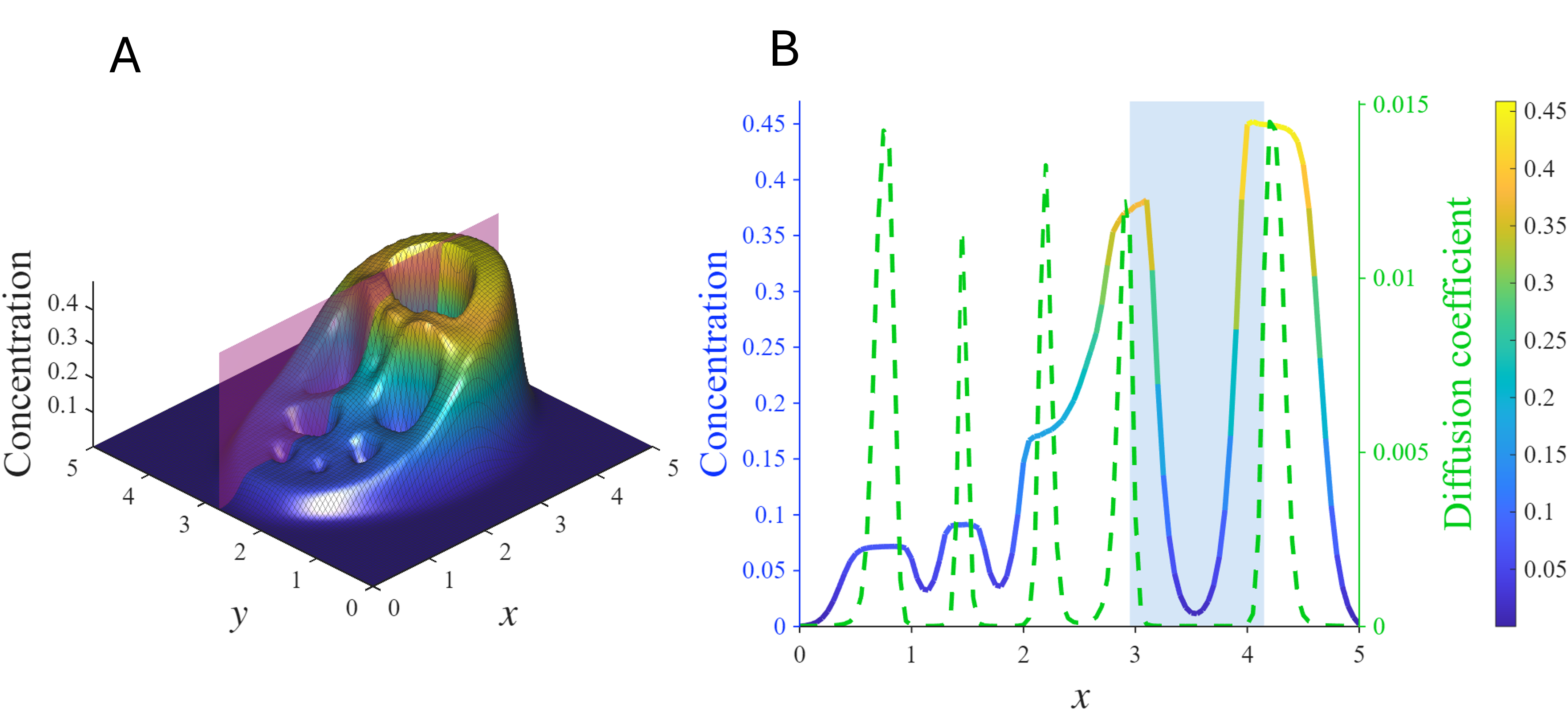}
\caption{\allblack 
{\bf Steady-state chemoattractant concentration and geometry-dependent
diffusion prior to border-cell migration.}
(A) Steady-state two-dimensional chemoattractant concentration field,
shown as a surface over the egg-chamber domain. Concentration is highest
near the oocyte at the posterior end and decreases toward the anterior (x=0),
while the cellular geometry produces spatial heterogeneity in the field.
The magenta plane marks the cross-section at $y=2.75$ analyzed in
panel~(B).
(B) Chemoattractant concentration and effective diffusion coefficient
along the same cross-section. The dashed green curve represents
$D(\boldsymbol{\phi})$, and the light-blue shaded region indicates the
spatial extent of the oocyte along this cross-section. Cellular regions correspond to strongly reduced
diffusivity, whereas extracellular gaps permit greater diffusion.
The resulting concentration profile therefore reflects both proximity
to the oocyte source and the heterogeneous diffusion landscape generated
by the cellular geometry (See Supplementary videos S2 and S4).
}
    \label{fig:cross sectional}
\end{figure}

{\allblack Figure~\ref{fig:cross sectional} illustrates the relationship between this geometry-dependent diffusivity and the resulting steady-state chemoattractant field. Panel~(A) shows the two-dimensional concentration field represented as a surface over the egg-chamber domain. The chemoattractant concentration is highest near the oocyte at the posterior end and decreases toward the anterior. However, this decrease is not spatially uniform: the cellular geometry produces channels and barriers to diffusion, resulting in a heterogeneous concentration landscape. The magenta plane identifies the cross-section at $y=2.75$ used in panel~(B).}
% \begin{figure}[htp!]
%     \centering
%     \includegraphics[width=1\linewidth]{Figs/Figure4_main.png}
% \caption{\allblack 
% {\bf Steady-state chemoattractant concentration and geometry-dependent
% diffusion prior to border-cell migration.}
% (A) Steady-state two-dimensional chemoattractant concentration field,
% shown as a surface over the egg-chamber domain. Concentration is highest
% near the oocyte at the posterior end and decreases toward the anterior (x=0),
% while the cellular geometry produces spatial heterogeneity in the field.
% The magenta plane marks the cross-section at $y=2.75$ analyzed in
% panel~(B).
% (B) Chemoattractant concentration and effective diffusion coefficient
% along the same cross-section. The dashed green curve represents
% $D(\boldsymbol{\phi})$, and the light-blue shaded region indicates the
% spatial extent of the oocyte along this cross-section. Cellular regions correspond to strongly reduced
% diffusivity, whereas extracellular gaps permit greater diffusion.
% The resulting concentration profile therefore reflects both proximity
% to the oocyte source and the heterogeneous diffusion landscape generated
% by the cellular geometry (See Supplementary videos S2 and S4).
% }
%     \label{fig:cross sectional}
% \end{figure}
{\allblack Panel~(B) provides a direct comparison between the chemoattractant
concentration and the effective diffusion coefficient along this
cross-section. The dashed green curve shows $D(\boldsymbol{\phi})$, while the concentration profile is shown on the left axis. Regions occupied by cells correspond to strongly reduced diffusivity, whereas extracellular gaps correspond to larger values of the diffusion coefficient. The concentration profile therefore reflects not only distance from the oocyte source but also the spatial variation in diffusivity generated by the cellular arrangement. In particular, changes in $D(\boldsymbol{\phi})$ along the cross-section are associated with localized changes in the concentration profile.
Thus, Fig.~\ref{fig:cross sectional} serves as a diagnostic of the
mechanism by which the phase-field geometry shapes chemoattractant
distribution. In contrast, the following figures examine the full spatial
distribution and temporal evolution of the concentration field during the migration simulations.
}

As an independent check of the diffuse-interface transport formulation, we
also compared the chemoattractant field with a sharp-interface finite-element
calculation in which cellular interiors are explicitly excluded from the
diffusion domain. The two formulations produce qualitatively similar
large-scale concentration patterns and geometry-dependent spatial structure,
with differences concentrated primarily near the differently represented
cellular interfaces (Supplementary Fig.~S4).

% Building on the spatially constrained diffusion described earlier, we now examine the evolution of the chemoattractant concentration field \( c(x,y,t) \) over time. 

% Figure~\ref{fig:cross sectional} (A) shows the one dimensional cross-sectional of concentration distribution (Eq.~\eqref{eq: chemo for oocyte}) and spatial diffusion (Eq.~\eqref{eq pfm diff}).  
% There is no diffusion of chemoattractant concentration into the cells (nurse cells and oocyte). Figure~\ref{fig:cross sectional}(B) represents the two dimensional chemoattractant concentration distribution corresponding to the steady state. The concentration is higher closer to the oocyte (posterior) and decreases to the anterior. Due to the heterogeneous diffusion landscape shaped by cell geometry, the resulting concentration field develops in a nonuniform and anisotropic manner.

\begin{figure}[htp!]
    \centering
    \includegraphics[width=1\linewidth]{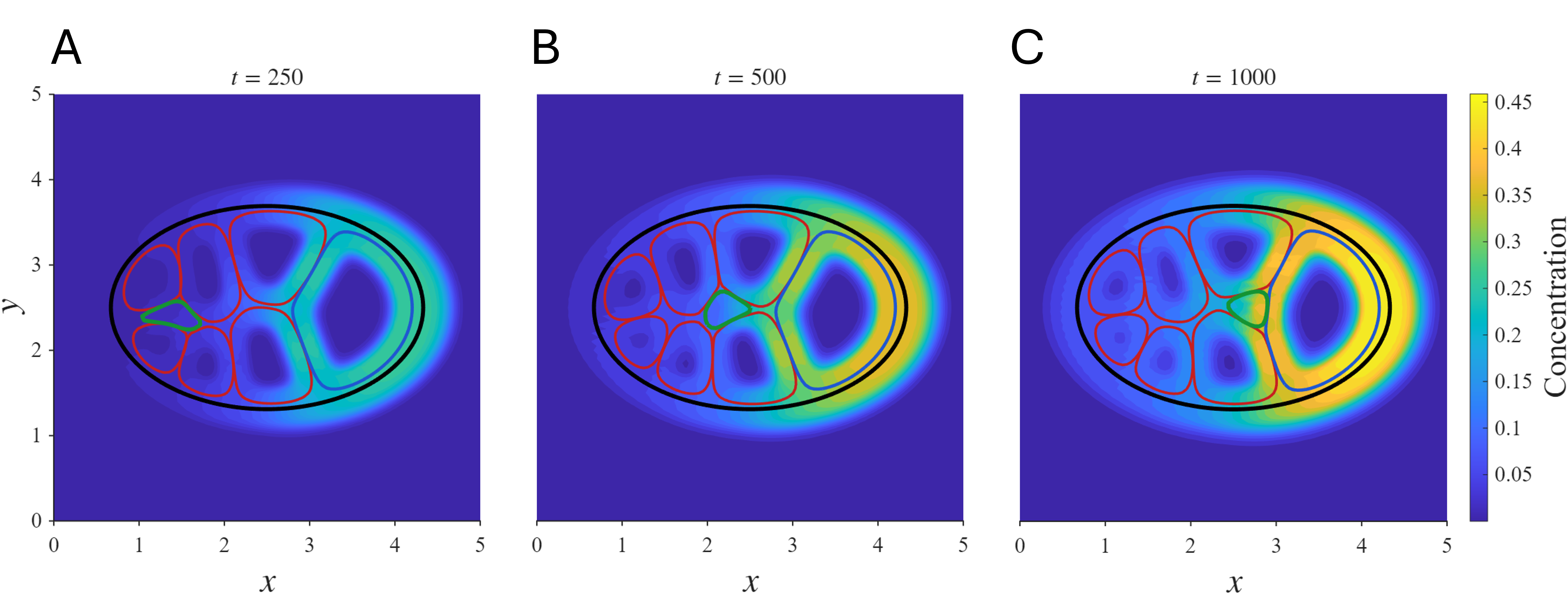}
    %  \caption{{\bf{Time evolution of the chemoattractant concentration field \( c(x,y,t) \)}}. Snapshots of the concentration field are shown $t= 115$ (A), $t=480$ (B), and $t=1000$ (C). The source of chemoattractant is localized to the posterior, corresponding to the oocyte boundary (right side), and diffusion occurs through the extracellular space while being excluded from the interiors of nurse cells and the oocyte. Over time, the distribution broadens and develops into a posterior–anterior gradient. By $t=1000$, the concentration field has stabilized into a steady-state profile, with the highest levels adjacent to the oocyte and progressively lower concentrations toward the anterior. The snapshots show how tissue geometry shapes anisotropic diffusion, producing heterogeneous gradients that provide directional cues for border cell migration (see in  Supplementary S4 Video).
    % }
    \caption{\allblack 
{\bf Time evolution of the chemoattractant concentration field for the dynamic case initialized from $c_0=0$.}
Snapshots of $c(x,y,t)$ are shown at (A) $t=250$, (B) $t=500$, and (C) $t=1000$. The $\phi=0.5$ level sets of the cellular phase fields are superimposed on the concentration field: nurse-cell boundaries are shown in red, the border-cell-cluster boundary in green, the oocyte boundary in blue, and the epithelial boundary in black.
The concentration field is displayed over the full computational grid for visualization; values outside the black epithelial contour are not interpreted
as part of the modeled egg-chamber concentration domain.
Chemoattractant is produced along the oocyte interface and spreads from the posterior toward the anterior while diffusion is strongly reduced within cellular regions and occurs preferentially through the extracellular spaces between cells. As time progresses, the concentration field broadens and develops a pronounced posterior-to-anterior gradient, while local variations in the field reflect the heterogeneous and evolving cellular geometry. The snapshots illustrate the concurrent evolution of the chemoattractant distribution and the phase-field geometry during border-cell migration (See Supplementary videos S2 and S4).
}
    \label{fig:concentration no intension}
\end{figure}

{\allblack
Figure~\ref{fig:concentration no intension} shows the temporal evolution of the dynamically evolving chemoattractant field for the case initialized from $c_0=0$. The concentration field is shown together with the $\phi=0.5$ level sets of the phase-field variables, which provide an approximate representation of the cellular boundaries. Nurse-cell boundaries are shown in red, the border-cell-cluster boundary in green, the oocyte boundary in blue, and the epithelial boundary in black. Thus, each panel shows both the evolving cellular geometry and the corresponding chemoattractant distribution.
{\allblack
Although the concentration field is plotted over the full rectangular computational grid, our analysis is restricted to the egg-chamber region enclosed by the epithelial boundary; concentration values displayed outside this contour are not used in the biological interpretation.
}

At the early stage, $t=250$, chemoattractant produced along the oocyte
interface has begun to diffuse into the surrounding extracellular domain according to Eq.~\eqref{eq: chemo for oocyte}. The concentration remains largest near the posterior oocyte boundary and decreases toward the anterior. Because the effective diffusion coefficient is strongly reduced within cellular regions, transport occurs preferentially through the extracellular spaces between the oocyte, nurse cells, border-cell cluster, and epithelial boundary, producing a spatially heterogeneous concentration field.

By $t=500$, the chemoattractant has propagated farther toward the anterior while the border-cell cluster has migrated into the nurse-cell environment. The concentration field therefore evolves in a geometry that is itself changing with time. The cellular boundaries restrict the available diffusion pathways, leading to nonuniform concentration contours and locally altered concentration levels around the intervening cells and extracellular gaps.

At the later time, $t=1000$, the chemoattractant field extends throughout a larger portion of the egg chamber and exhibits a pronounced posterior-to-anterior concentration gradient, with the largest concentrations remaining near the oocyte source. The persistence of lower concentrations within cellular interiors and higher concentrations along extracellular regions reflects the phase-field-dependent diffusivity introduced in Eq.~\eqref{eq pfm diff}. These snapshots therefore illustrate how the time-dependent chemoattractant field develops from the initially zero-concentration state while being continuously shaped by the evolving cellular geometry.
}

\subsection{Directed Migration of the Border Cell Cluster}
In classical models of chemotaxis, the directed motion of a cell or cluster is described by a chemical force \cite{nonomura2012study}.
% of the form
% \begin{align*}
%     \mathbf{F}_{\text{chem}} = -\mu_c \nabla \cdot (\phi_c \nabla c)
% \end{align*}
% where $\mu_c$ is the chemotactic sensitivity, $\phi_c$ is the phase-field variable representing the cluster, and $c$ is the chemoattractant concentration. 
This force biases movement along concentration gradients and has been widely applied in both single- and multi-cell models \cite{alt1980biased, horstmann20031970, hillen2009user}. However, such a force assumes that cells respond uniformly to extracellular signals across their entire boundary and does not typically capture the critical role of physical interactions with neighboring substrate cells in the region.

To incorporate these biological features, we instead employ the Tangential Interface Migration (TIM) force (Eq.~\eqref{eq: tim force corrected}), which acts specifically at the interfaces where the cluster contacts adjacent nurse cells. The TIM force couples local chemoattractant gradients with interfacial geometry along cell–cell contacts rather than uniformly across the cluster boundary. 
Compared to the traditional chemical force, TIM therefore provides a more realistic representation of mechanically constrained, contact-mediated migration \emph{in vivo}.
% {\allcpink (i think the two paragraphs above could come before this section when you first talk about the TIM model earlier. but up to you)}{\allorange Does feel like it should go prior to introducing TIM}

\begin{figure}[htp!]
    \centering
    \includegraphics[width=1\linewidth]{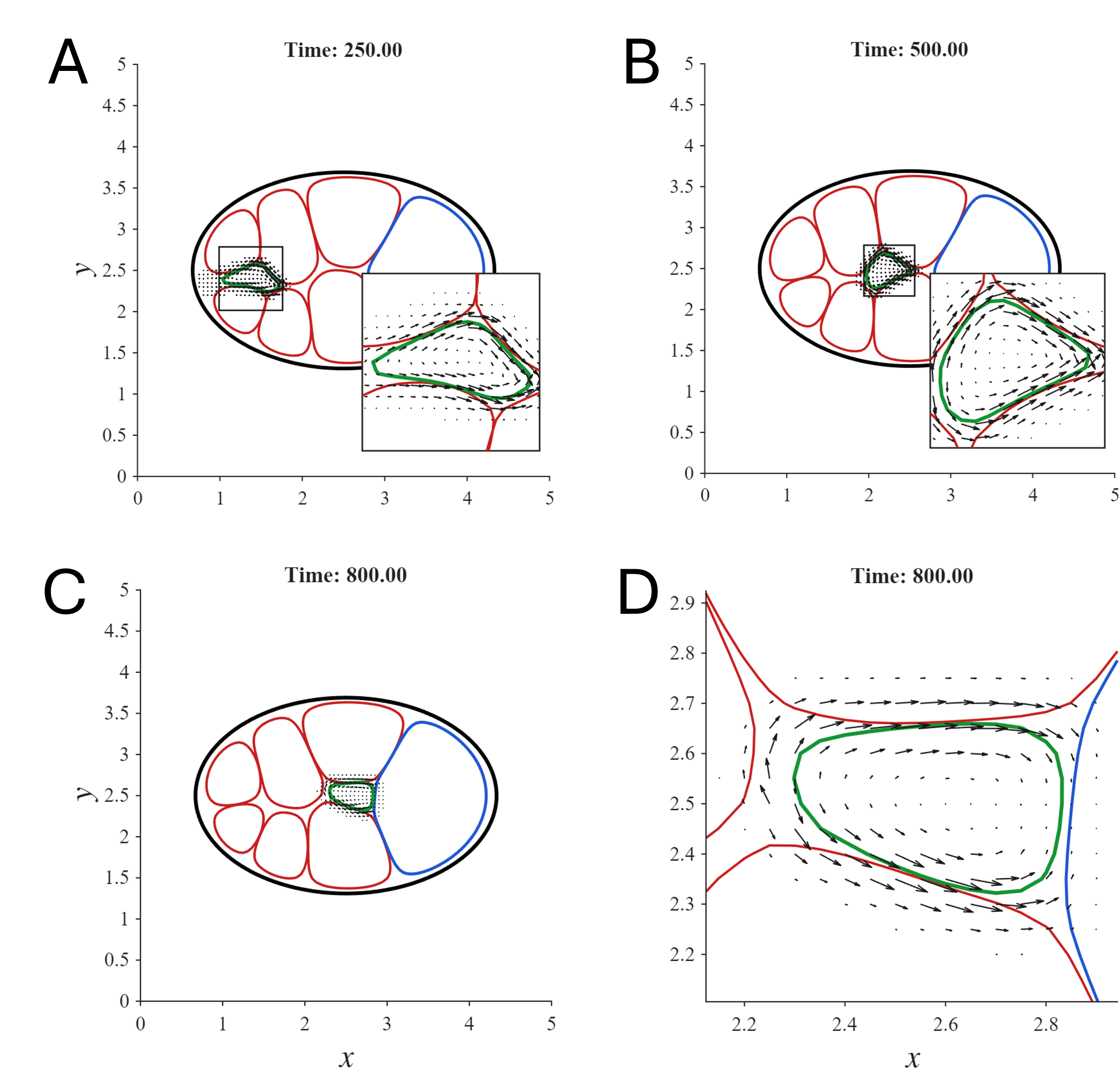}
%         \caption{{\bf{Time evolution of border cell cluster migration in the {Drosophila} egg chamber simulated using the phase field model.}} Cell boundaries are represented by phase-field contours (red, nurse cells; blue, oocyte; black, epithelial layer). The border cell cluster (green) initiates at the anterior end and migrates toward the posterior oocyte under the influence of the tangential interface migration (TIM) force. The model incorporates adhesive interactions to maintain cluster cohesion, repulsive forces to prevent overlap, and volume constraints to preserve cell size. Magenta arrows indicate the magnitude and direction of the TIM force.
% (A) At $t=115$, the cluster begins to respond to the TIM force and initiates migration after the chemoattractant concentration reaches steady state.
% (B) At $t=310$, the TIM force acts in a contact-mediated manner at interfaces between the cluster and adjacent nurse cells, driving forward movement.
% (C) At $t=785$, the cluster has nearly completed migration and establishes contact with the oocyte. (D) The enlarged panel of the cluster migration at $t=785$ shows the tangentially oriented vectors along the cluster boundary (See in  Supplementary S5 Video).}
\caption{\allblack 
{\bf Time evolution of border-cell-cluster migration and TIM forcing in the
\emph{Drosophila} egg chamber.}
Cellular boundaries are represented by the $\phi=0.5$ level sets of the
phase-field variables: nurse-cell boundaries are shown in red, the
border-cell-cluster boundary in green, the oocyte boundary in blue, and the
epithelial boundary in black. Black arrows show the local tangential
interface migration (TIM) forcing acting in the vicinity of the
border-cell-cluster interface.
(A) At $t=250$, the cluster has begun migrating from the anterior region
through the surrounding nurse-cell environment, with the TIM vectors
localized around the cluster boundary.
(B) At $t=500$, the cluster has progressed farther toward the posterior, and
the TIM forcing changes orientation as the cluster interacts with neighboring
nurse-cell interfaces. The inset provides an enlarged view of the cluster and
the local force-vector organization.
(C) At $t=800$, the cluster has approached the oocyte and occupies the
posterior nurse-cell region.
(D) Enlarged view of the cluster at $t=800$, showing the local TIM vectors
along and around the cluster interface as it interacts with the surrounding
nurse-cell and oocyte boundaries.
Together, the snapshots illustrate how the TIM forcing evolves with the
changing cluster geometry and local cell--cell configuration during migration
(see Supplementary Video~S5).
}
        \label{fig:contour migration with tim}
\end{figure}

% Figure~\ref{fig:contour migration with tim} illustrates the time evolution of the border cell cluster under the combined influence of adhesive interactions, repulsive forces, volume constraints, and the TIM force. At $t=115$, the cluster begins to respond to the TIM force once the chemo- attractant concentration has reached steady state. By $t=310$, the TIM force acts at contact interfaces between the cluster and surrounding nurse cells, driving forward motion in a directed manner. At $t=785$, the cluster has nearly completed migration and establishes contact with the oocyte. Together, these snapshots highlight how TIM integrates chemical guidance and mechanical constraints to produce persistent, directed migration.
{\allblack Figure~\ref{fig:contour migration with tim} illustrates the time evolution of border-cell-cluster migration under the combined influence of adhesive interactions, repulsive forces, volume constraints, and the TIM force. At $t=250$, the cluster has begun migrating from the anterior region through
the surrounding nurse-cell environment, with the TIM forcing localized around the cluster interface. By $t=500$, the cluster has progressed farther toward the posterior, and the orientation of the TIM vectors changes as the cluster interacts with neighboring nurse-cell interfaces. At $t=800$, the cluster has approached the oocyte and occupies the posterior nurse-cell region. The enlarged views in Fig.~\ref{fig:contour migration with tim} highlight the local organization of the TIM vectors along the cluster boundary at different stages of migration. Together, these snapshots illustrate how the contact-mediated TIM forcing evolves with the changing cluster geometry and local cell--cell configuration while supporting directed migration toward the oocyte.}
{\allblack The simulated migration is qualitatively consistent with live-imaging observations in which the border cell cluster moves posteriorly through the nurse-cell environment while maintaining extensive contacts with surrounding cells and ultimately reaches the oocyte \cite{dai2020tissue,akhavan2026phase}. The present simulations are not fitted quantitatively to individual experimental trajectories or cell shapes; rather, they are intended to reproduce the principal geometry- and contact-dependent features of border-cell migration while examining how these features interact with the spatial chemoattractant field.
}

{\allblack
As an additional check of the diffuse-interface formulation, we compared
the resulting migration dynamics with those obtained from the corresponding
sharp-interface finite-element model. Both formulations reproduce the same
overall posterior progression of the border-cell cluster and qualitatively
similar migration dynamics, although quantitative differences remain in the
instantaneous speed and final cluster position. This comparison indicates
that the principal migration behavior is retained when the diffuse
phase-field representation is replaced by explicitly resolved sharp cellular
boundaries (Supplementary Fig.~S5).
}

\subsection{Role of Chemoattractant Field Dynamics in Border Cell Migration}
\begin{figure}[htp!]
    \centering
    \includegraphics[width=\linewidth]{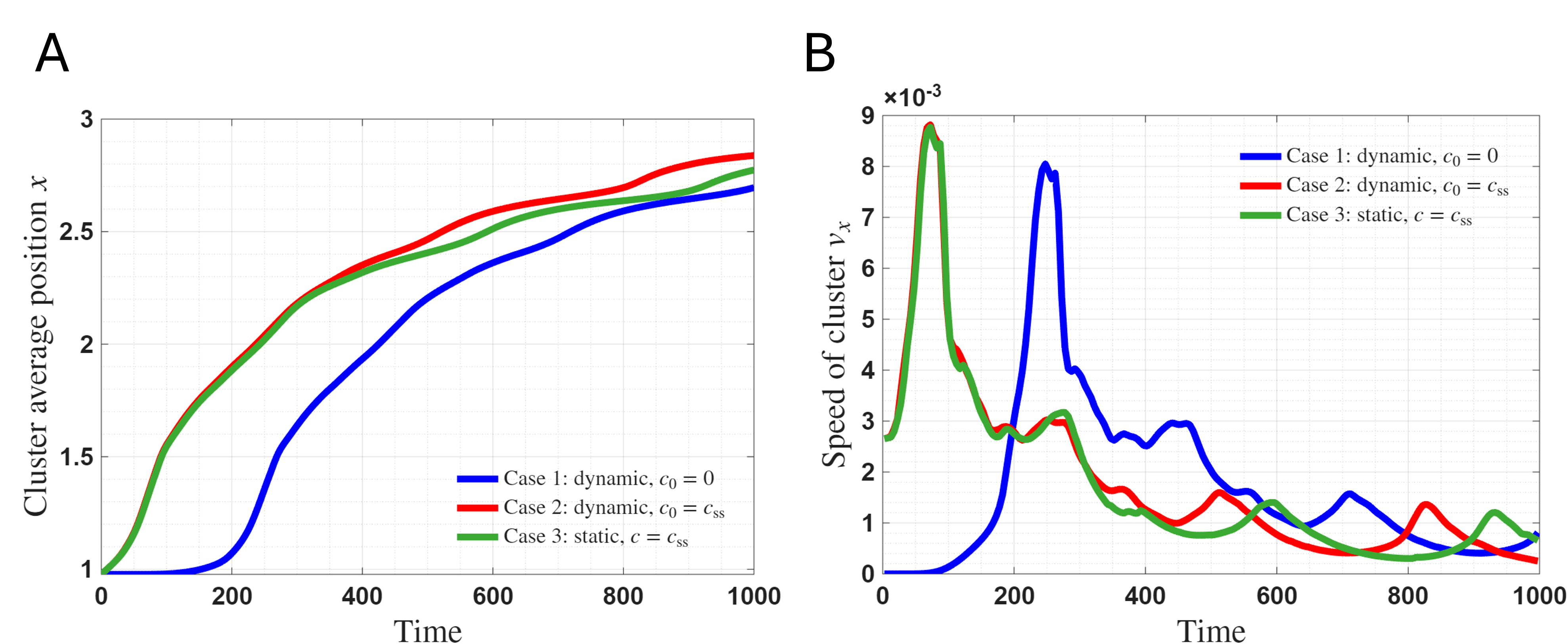}
    \caption{
{\bf Comparison of migration under dynamic and static chemoattractant fields.}
Three concentration conditions are considered: Case~1, a dynamically
evolving concentration field initialized at $c_0=0$; Case~2, a
dynamically evolving field initialized from the steady-state
concentration $c_0=c_{\mathrm{ss}}$; and Case~3, a static concentration
field fixed at $c=c_{\mathrm{ss}}$. (A) Average $x$-position of the
border-cell cluster as a function of time. (B) Signed migration speed
$v_x$ as a function of time. The baseline dynamic case exhibits a clear
delay while the chemoattractant field develops, whereas the dynamic and
static cases initialized from $c_{\mathrm{ss}}$ show substantially more
similar migration behavior.
}
\label{fig:dynamic_static_comparison}
\end{figure}

{\allblack
To determine how the temporal evolution of the chemoattractant field influences border-cell-cluster migration, we compared three concentration conditions while keeping the phase-field geometry and migration parameters unchanged. We denote by $c_{\mathrm{ss}}(\mathbf{x})$ the steady-state chemoattractant field obtained from the diffusion--reaction model prior to cluster migration.

In the first case, corresponding to the baseline formulation, the chemoattractant concentration is initialized at $c_0=0$ and evolves dynamically throughout the simulation. In the second case, the concentration is initialized from the steady-state field, $c_0=c_{\mathrm{ss}}$, and subsequently allowed to evolve dynamically as the cluster moves. In the third case, the same steady-state field is held fixed throughout migration, such that $c(\mathbf{x},t)=c_{\mathrm{ss}}(\mathbf{x})$.

To quantify the transient establishment of the chemoattractant field, we
also compared the solution initialized from $c_0=0$ with an independently
computed steady-state solution $c_{\mathrm{ss}}$. The relative $L^2$ and
$L^\infty$ errors decrease progressively with time and reach the percent
level only at approximately $t\approx2000$ (Supplementary Fig.~S6).
Accordingly, the concentration field is still far from steady state when
cluster migration begins and continues to evolve substantially during the
early migration period. Cross-sectional concentration profiles at selected
times further illustrate this convergence toward $c_{\mathrm{ss}}$.

Figure~\ref{fig:dynamic_static_comparison}(A) compares the average $x$-position of the cluster for the three cases. Initializing the concentration from $c_0=0$ produces a pronounced delay in migration relative to the two cases initialized from $c_{\mathrm{ss}}$. By contrast, the dynamic and static simulations initialized from the steady-state field follow substantially more similar trajectories, although a modest separation develops as migration proceeds.

The corresponding migration speeds are shown in Fig.~\ref{fig:dynamic_static_comparison}(B). The largest differences occur during the early transient regime. In the baseline dynamic case, cluster acceleration is delayed until the chemoattractant field develops, whereas the two steady-state-initialized cases exhibit an earlier increase in migration speed. At later times, the dynamic and static steady-state-initialized cases display broadly similar speed profiles with smaller local differences.

% \begin{figure}[htp!]
%     \centering
%     \includegraphics[width=\linewidth]{Figs/Figure7_main.png}
%     \caption{
% {\bf Comparison of migration under dynamic and static chemoattractant fields.}
% Three concentration conditions are considered: Case~1, a dynamically
% evolving concentration field initialized at $c_0=0$; Case~2, a
% dynamically evolving field initialized from the steady-state
% concentration $c_0=c_{\mathrm{ss}}$; and Case~3, a static concentration
% field fixed at $c=c_{\mathrm{ss}}$. (A) Average $x$-position of the
% border-cell cluster as a function of time. (B) Signed migration speed
% $v_x$ as a function of time. The baseline dynamic case exhibits a clear
% delay while the chemoattractant field develops, whereas the dynamic and
% static cases initialized from $c_{\mathrm{ss}}$ show substantially more
% similar migration behavior.
% }
% \label{fig:dynamic_static_comparison}
% \end{figure}

These comparisons separate two effects of the chemoattractant model. First, the transient establishment of the guidance field can strongly affect the onset and timing of migration. Second, once the chemoattractant distribution is already close to steady state, continuously updating the concentration field produces a more modest effect relative to maintaining the same geometry-shaped field as a static profile. Thus, in the parameter regime considered here, the principal consequence of the time-dependent diffusion--reaction model arises during formation of the guidance field, while a static steady-state profile provides a reasonable approximation for much of the subsequent migration.

% These comparisons separate two effects of the chemoattractant model. First,
% the transient establishment of the guidance field can strongly affect the
% onset and timing of migration. Second, once the chemoattractant distribution
% is already close to steady state, continuously updating the concentration
% field produces a more modest effect relative to maintaining the same
% geometry-shaped field as a static profile. Thus, in the parameter regime
% considered here, the principal consequence of the time-dependent
% diffusion--reaction model arises during formation of the guidance field,
% while a static steady-state profile provides a reasonable approximation for
% much of the subsequent migration.

The migration speed nevertheless remains strongly nonuniform as the cluster
moves through the egg chamber. To examine whether these local changes are
associated with the surrounding tissue geometry, we compared selected local
maxima and minima of $v_x(t)$ with the corresponding cellular configurations
and concentration fields (Supplementary Fig. S2). Local speed minima,
including those near $t=400$ and $t=640$, coincide with more confined
configurations in which the cluster is constrained between neighboring
nurse-cell interfaces, whereas local increases in speed occur as the cluster
moves through or emerges from these constricted regions. This correspondence
shows that the heterogeneous migration speed is closely associated with the
sequence of geometric bottlenecks encountered along the migration path.
Because these geometric changes affect both mechanical confinement and the
phase-dependent diffusion landscape, we do not attribute each individual
speed variation exclusively to the chemoattractant field. The local coupling
between interface motion and concentration dynamics is examined separately
in the following subsection.
}

\subsection{Local Feedback Between Moving Interfaces and the Chemoattractant Field}

\begin{figure}[htp!]
    \centering
    \includegraphics[width=1\linewidth]{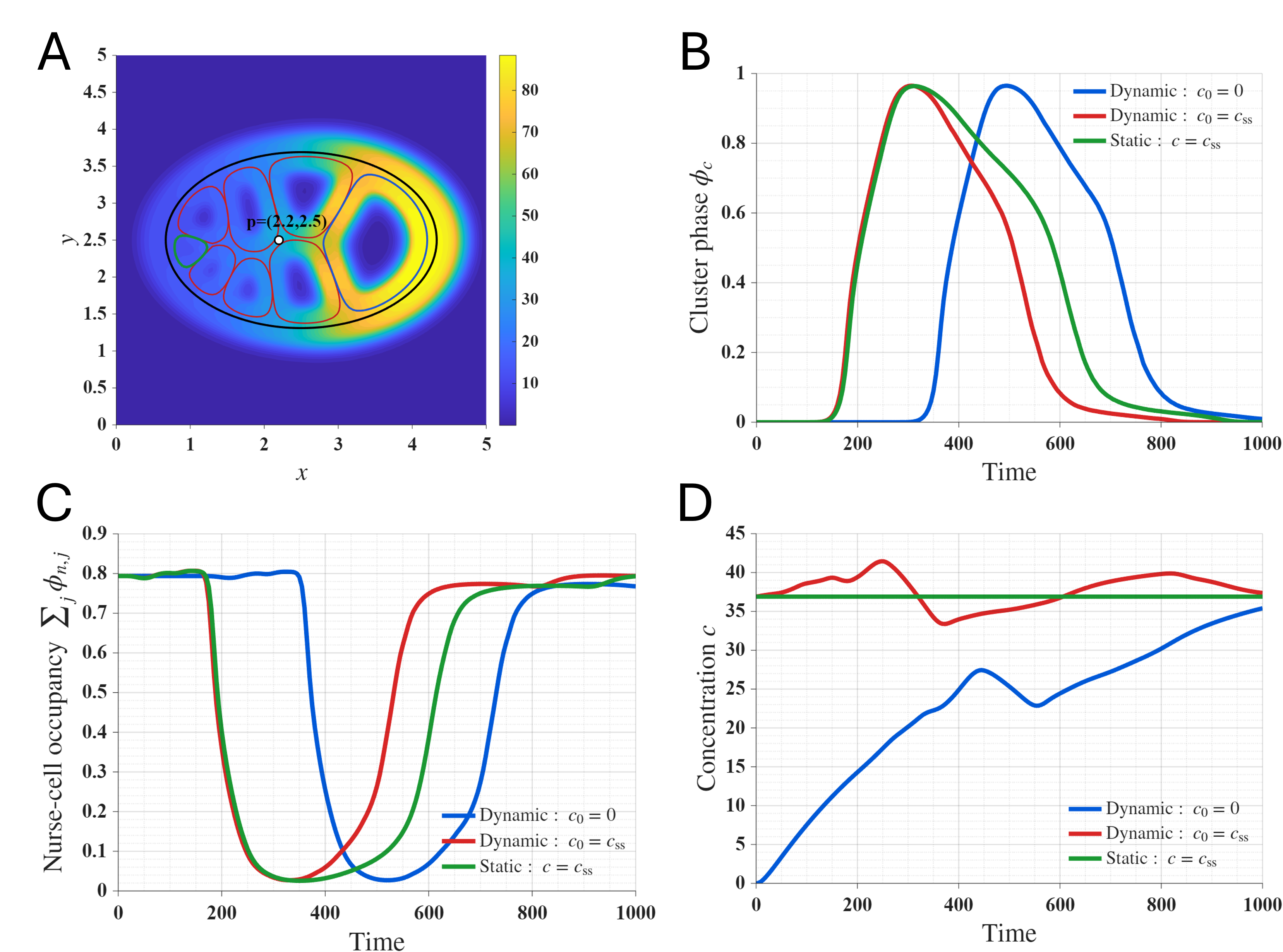}
    \caption{
{\bf Local coupling between moving cellular interfaces and the
chemoattractant field.}
All quantities are evaluated at the fixed probe point
$p=(2.2,2.5)$ for the three concentration conditions introduced in the
preceding subsection: Case~1, a dynamically evolving concentration field
initialized from $c_0=0$; Case~2, a dynamically evolving field initialized
from $c_0=c_{\mathrm{ss}}$; and Case~3, a static field fixed at
$c=c_{\mathrm{ss}}$.
(A) Local border-cell-cluster phase
$\phi_{\mathrm{cluster}}(p,t)$, identifying the approach, passage, and
departure of the cluster through the probe region.
(B) Summed nurse-cell phase
$\sum_i\phi_{{\rm nurse},i}(p,t)$, showing the corresponding local
rearrangement of the surrounding nurse-cell interfaces.
(C) Local chemoattractant concentration $c(p,t)$. In the static case the
concentration remains fixed, whereas in the dynamically evolving cases it
changes as the cellular geometry evolves. The temporal correspondence between
the phase-field changes and concentration variations demonstrates the local
coupling between moving interfaces and geometry-dependent chemoattractant
transport.
}
\label{fig:local_feedback}
\end{figure}

{\allblack
The preceding comparison showed that the temporal treatment of the
chemoattractant field influences the onset and subsequent trajectory of
border-cell-cluster migration. To examine how these differences arise
locally, we next considered the coupling between the moving cellular
interfaces and the chemoattractant field at a fixed spatial location.
Specifically, we monitored the cluster phase field, the summed nurse-cell
phase fields, and the chemoattractant concentration at the probe point
$p=(2.2,2.5)$ for the same three concentration conditions considered in
the previous subsection.

% \begin{figure}[htp!]
%     \centering
%     \includegraphics[width=1\linewidth]{Figs/Figure8_main.png}
%     \caption{
% {\bf Local coupling between moving cellular interfaces and the
% chemoattractant field.}
% All quantities are evaluated at the fixed probe point
% $p=(2.2,2.5)$ for the three concentration conditions introduced in the
% preceding subsection: Case~1, a dynamically evolving concentration field
% initialized from $c_0=0$; Case~2, a dynamically evolving field initialized
% from $c_0=c_{\mathrm{ss}}$; and Case~3, a static field fixed at
% $c=c_{\mathrm{ss}}$.
% (A) Local border-cell-cluster phase
% $\phi_{\mathrm{cluster}}(p,t)$, identifying the approach, passage, and
% departure of the cluster through the probe region.
% (B) Summed nurse-cell phase
% $\sum_i\phi_{{\rm nurse},i}(p,t)$, showing the corresponding local
% rearrangement of the surrounding nurse-cell interfaces.
% (C) Local chemoattractant concentration $c(p,t)$. In the static case the
% concentration remains fixed, whereas in the dynamically evolving cases it
% changes as the cellular geometry evolves. The temporal correspondence between
% the phase-field changes and concentration variations demonstrates the local
% coupling between moving interfaces and geometry-dependent chemoattractant
% transport.
% }
% \label{fig:local_feedback}
% \end{figure}

Figure~\ref{fig:local_feedback}(A) shows the local cluster phase
$\phi_{\mathrm{cluster}}(p,t)$. For the simulations initialized from the
steady-state concentration field, the cluster reaches the probe location
substantially earlier than in the simulation initialized from $c_0=0$.
This is reflected by the earlier increase of
$\phi_{\mathrm{cluster}}$ in Cases~2 and~3. In Case~1, the corresponding
increase occurs later, consistent with the delayed migration observed in
the cluster-position and speed curves of the preceding subsection.
The subsequent decay of $\phi_{\mathrm{cluster}}$ indicates passage of
the cluster beyond the probe location.

The motion of the cluster is accompanied by a corresponding rearrangement
of the surrounding nurse-cell interfaces. As shown in
Fig.~\ref{fig:local_feedback}(B), the summed local nurse-cell phase
$\sum_i \phi_{{\rm nurse},i}(p,t)$ decreases as the cluster approaches
and occupies the probe region, and subsequently increases again after the
cluster passes. The timing of this decrease is shifted between the three
cases in a manner consistent with the local cluster-phase signal. Together,
panels~(A) and~(B) therefore provide a local measure of the evolving
cellular geometry around the probe point.

The effect of this interface motion on chemoattractant distribution is shown
in Fig.~\ref{fig:local_feedback}(C). In the static case,
$c(\mathbf{x},t)=c_{\mathrm{ss}}(\mathbf{x})$, the concentration at the
probe remains constant by construction. In contrast, when the concentration
is initialized from the same steady-state field but subsequently allowed to
evolve dynamically, the local concentration departs appreciably from its
initial value as the cellular configuration changes. The concentration first
increases above the steady-state value and later decreases before recovering,
with these variations occurring over the same interval in which the cluster
and nurse-cell phase fields undergo substantial changes. Thus, even when the
initial chemoattractant distribution is identical, motion of the cellular
interfaces produces a time-dependent redistribution of the local
chemoattractant field through the phase-dependent diffusivity.

The case initialized from $c_0=0$ exhibits the same coupling on a delayed
time scale. The chemoattractant field first develops from zero, while the
cluster reaches the probe later; changes in the local cluster and nurse-cell
phase fields are then accompanied by corresponding variations in the
concentration profile. These observations are consistent with the global
migration differences reported in the preceding subsection and demonstrate
the two-way coupling represented by the model: cellular-interface motion
modifies the geometry through which chemoattractant is transported, while
the resulting concentration field contributes to the migration force acting
on the border-cell cluster.

The magnitude of this feedback is sufficient to produce measurable local
changes in chemoattractant concentration and quantitative differences in
migration timing and trajectory, although the parameter regime considered
here does not lead to a loss of directed migration. Having established this
local coupling between cellular geometry and chemoattractant transport, we
next examine how changing a geometric property of the migrating group itself,
namely the target cluster volume, influences both migration dynamics and the
chemoattractant field.
}

\subsection{Influence of Cluster Volume on Migration Dynamics}
Collective migration is not only guided by chemical cues but also shaped by the physical size of the migrating group \cite{stonko2015mathematical, aranjuez2016dynamic, george2025chemotaxis, cai2016modeling}.
% {\allcpink (need to cite Cai and Stonko here)} {\allblack (Cai modeling paper?)}
For the \emph{Drosophila} border cell cluster, cluster size could influence migration in two competing ways.
% On one hand, larger clusters may provide an advantage in gradient sensing: the difference in chemoattractant concentration between the front and back increases with distance, potentially enhancing the cluster’s ability to polarize and maintain directionality \cite{george2025chemotaxis}.
% On the other hand, larger clusters displace more surrounding tissue and experience stronger mechanical resistance from nurse cells and the extracellular matrix, which act as both the substrate and the confining environment for migration \cite{cai2016modeling,aranjuez2016dynamic,lamb2021fascin, dai2020tissue}. These competing effects suggest the existence of an optimal cluster size that balances chemical sensitivity with physical constraints.
{\allblack {\allblack
Cluster size can influence migration through both mechanical and chemical
mechanisms. Mechanically, larger clusters displace more surrounding tissue
and can experience greater resistance from neighboring nurse cells and the surrounding
extracellular matrix, which together form the substrate and confining
environment for migration
\cite{cai2016modeling,aranjuez2016dynamic,dai2020tissue, stonko2015mathematical, molina2023constriction}.}
% {\allcpink{Molina López E, Kabanova A, Winkel A, Franze K, Palacios IM, Martín-Bermudo MD. Constriction imposed by basement membrane regulates developmental cell migration. PLoS Biol. 2023 Jun 28;21(6):e3002172. doi: 10.1371/journal.pbio.3002172. PMID: 37379333; PMCID: PMC10335704.}
Chemically, because the border-cell cluster itself contributes to the
geometry-dependent diffusion landscape, changing its size can alter the
spatial distribution of chemoattractant in the surrounding migration
domain. Thus, rather than assuming that cluster size affects migration
solely through mechanical confinement or through an intrinsic advantage
in gradient sensing, we examine how changes in cluster size modify both
migration dynamics and the local chemoattractant field.
}

% \begin{figure}[ht!]
%     \centering
%     \includegraphics[width=1\linewidth]{Figs/figure7.png}
%     \caption{{\bf{Effect of cluster volume on migration dynamics.}}
% (A) Position of the cluster center over time for three target volumes (0.12, 0.15, 0.18). The larger cluster (green) shows slower displacement once migration begins, indicating increased difficulty in squeezing between nurse cells compared to smaller clusters.
% (B) Instantaneous speed of the cluster. The smaller cluster (blue) reaches the highest peak velocity, while the larger cluster (green) consistently migrates more slowly, reflecting the greater mechanical resistance imposed by surrounding nurse cells. These results suggest that larger clusters experience more hindrance in confined tissue environments, reducing their overall migration efficiency relative to smaller clusters.}
%     \label{fig:diff volumes}
% \end{figure}
The influence of cluster size on migration has been explored by Cai et al., who combined live imaging with a 
% {\allcpink (maybe something other than "simple" here- agent based? minimal?)}{\allorange maybe reduced}
reduced theoretical model of border cell migration \cite{cai2016modeling}. Their analyses suggested that migration speed depends on cluster size in a non-monotonic manner, with moderate increases in size enhancing migration but over-sized clusters becoming hindered by drag and confinement.

To test how cluster size influences migration within our chemo-mechanical
phase-field model, we systematically varied the target volume of the
border-cell cluster while keeping all other parameters fixed. The
intermediate case, $\bar V_c=0.15$ (magenta), was taken as the baseline,
while $\bar V_c=0.12$ (blue) and $\bar V_c=0.18$ (green) represent
approximately $20\%$ decreases and increases in target cluster size,
respectively.

\begin{figure}[ht!]
    \centering
    \includegraphics[width=\linewidth]{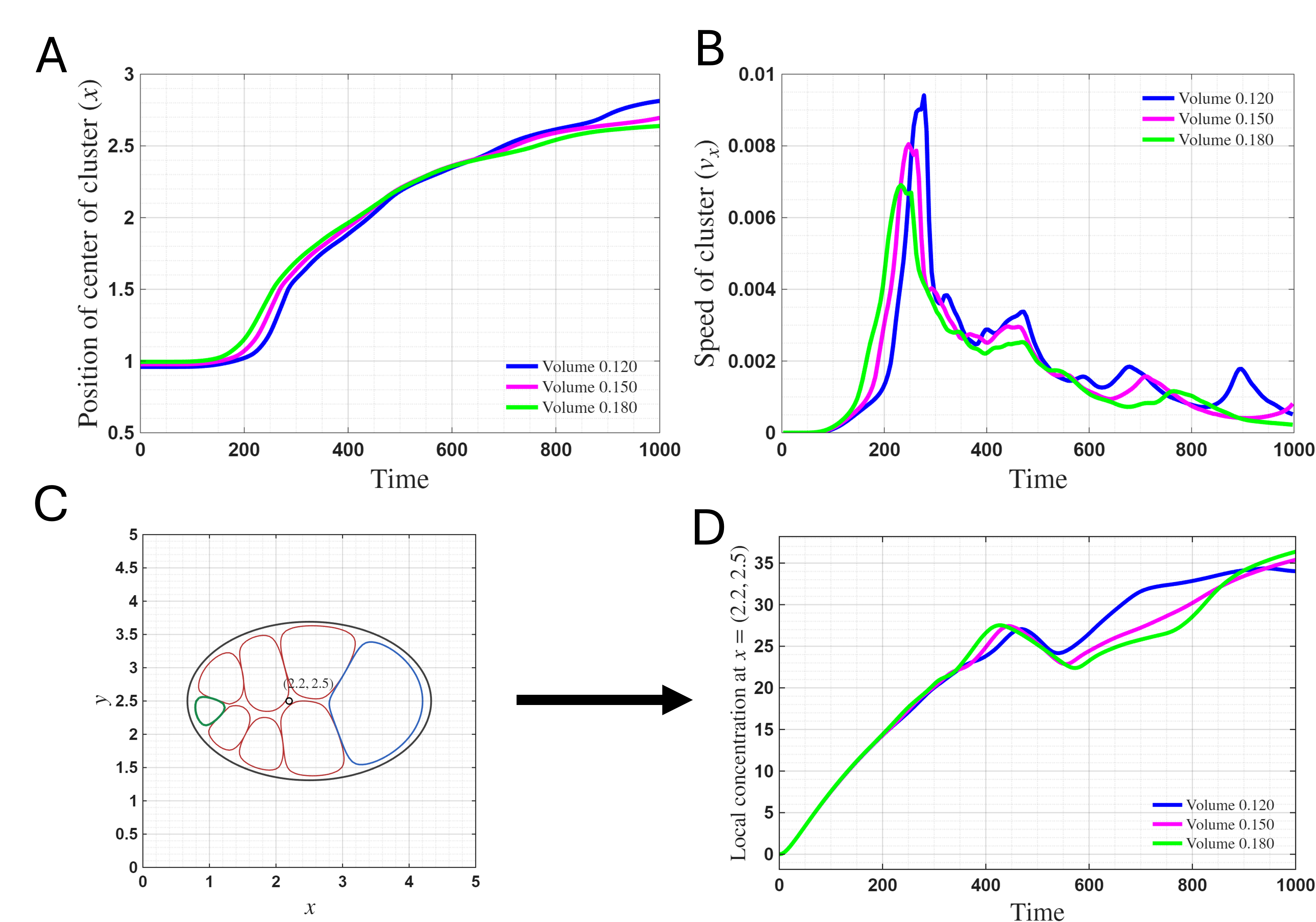}
    \caption{
{\bf Effect of cluster volume on migration dynamics and the local chemoattractant field.}
(A) Position of the cluster center over time for three target cluster volumes,
$\bar V_c=0.12$, $0.15$, and $0.18$. Although the larger cluster begins to
move slightly earlier, its subsequent displacement is reduced relative to the
smaller cluster, which reaches the largest final $x$-position.
(B) Instantaneous migration speed of the cluster. The smaller cluster
($\bar V_c=0.12$) reaches the highest peak speed, whereas increasing cluster
volume reduces the peak migration speed and shifts its timing, consistent with
greater mechanical hindrance during passage through the confined nurse-cell
environment.
(C) Egg chamber geometry showing the fixed probe location
$p=(2.2,2.5)$ used to evaluate the chemoattractant field. The cluster-center \(x\)-position reaches the probe location \(x=2.2\) at approximately \(t=487\) for \(\bar V_c=0.15\), \(t=500\) for \(\bar V_c=0.18\), and \(t=507\) for \(\bar V_c=0.12\).
(D) Local chemoattractant concentration $c(p,t)$ for the three cluster
volumes. The concentration profiles are nearly identical at early times but
separate as migration proceeds, and their ordering changes over time.
Thus, changing cluster volume affects not only cluster displacement and speed
but also the evolving local chemoattractant distribution through its effect on
the geometry-dependent diffusion field.
}
    \label{fig:diff volumes}
\end{figure}

{\allblack
Figure~\ref{fig:diff volumes} compares both the migration dynamics and the
local chemoattractant environment for these three cases. Panels~(A) and~(B)
show the cluster-center position and instantaneous migration speed,
respectively, while panels~(C) and~(D) show the probe location and the
corresponding local chemoattractant concentration.

The three cluster sizes produce distinct migration dynamics. The larger
cluster ($\bar V_c=0.18$) begins to move slightly earlier and reaches its
peak speed earlier than the baseline and smaller-cluster cases. However,
its peak speed is lower, and its subsequent displacement becomes smaller
than that of the smaller cluster. In contrast, the smaller cluster
($\bar V_c=0.12$) exhibits a later but larger peak in migration speed and
ultimately reaches the largest $x$-position by the end of the simulation.
The speed curves also show several later local variations and crossings,
indicating that the effect of cluster size is not simply a uniform
rescaling of migration speed throughout the trajectory. Rather, cluster
size modifies both the timing and magnitude of the cluster response as it
moves through the heterogeneous nurse-cell environment.

Because the border-cell cluster also contributes to the
geometry-dependent diffusion landscape, we next examined whether changing
cluster size modifies the chemoattractant field itself.
Figure~\ref{fig:diff volumes}(C) identifies the fixed probe point
$p=(2.2,2.5)$ within the migration region, and
Fig.~\ref{fig:diff volumes}(D) shows the corresponding concentration
$c(p,t)$ for the three cluster sizes. The concentration profiles are nearly
indistinguishable during the early stage, before the cluster configurations
around the probe have diverged substantially. As migration proceeds, the
profiles separate, with both the magnitude and ordering of the local
concentration depending on time. In particular, no single cluster size
produces consistently higher or lower concentration throughout the
simulation. This non-monotonic ordering indicates that cluster size affects
the local chemoattractant environment through the evolving cellular
geometry, rather than simply increasing or decreasing the concentration by
a fixed amount.
}

% The results show that the cluster size strongly influences migration dynamics, {i.e.}, larger clusters are mechanically hindered, while smaller clusters move more efficiently through the nurse cells environment. 
% % These results support the idea that there is an optimal cluster size for efficient migration: clusters that are too small may sense gradients less effectively, while overly large clusters are mechanically hindered. 
% Although our current simulations do not explore the lower size limit, where very small cluster might fail to maintain cohesion or directional sensing, the framework sets the stage for investigating how both excessively large and excessively small cluster sizes could impair migration.

% Very small clusters may also require stronger and more sustained adhesion–repulsion interactions with nurse cells and the oocyte to ensure effective movement.
% {\allblack These results demonstrate that cluster size can substantially modify the mechanical constraints experienced during migration and suggest that developmental regulation of cluster size may contribute to robust movement through the confined nurse-cell environment.}
{\allblack
Taken together, these results show that cluster size influences both the
mechanical and chemical components of migration. Larger clusters experience
greater mechanical hindrance while moving through the confined nurse-cell
environment, whereas changes in cluster size also modify the local
chemoattractant field through the geometry-dependent diffusion model. Thus,
the effect of cluster size in the coupled model cannot be attributed solely
to mechanical resistance; cluster geometry also feeds back on the chemical
environment available to the migrating cells. Our current simulations do not
explore substantially smaller clusters or the lower size regime in which
cohesion or directional sensing might become compromised, and we therefore
leave this regime for future investigation.
}

% Our phase-field framework captures this balance, reinforcing the conclusion that cluster size is a critical determinant of collective cell migration dynamics in the egg chamber.

{\allblack The dependence of migration on cluster size also raises the question of how the physical properties of the migrating group are regulated {\it in vivo}. Border-cell number is developmentally constrained through the signaling network that specifies migratory follicle-cell fate, particularly the JAK/STAT pathway and its positive and negative regulators \cite{monahan2013socs36e,montell2012group, hombria2002fertile}}.
% {\allcpink {Hombría JC, Brown S. The fertile field of Drosophila Jak/STAT signalling. Curr Biol. 2002 Aug 20;12(16):R569-75. doi: 10.1016/s0960-9822(02)01057-6. PMID: 12194841.}
Thus, the relatively restricted size of the normal border-cell cluster is not arbitrary. Previous experimental and modeling work has also demonstrated a non-monotonic relationship between border-cell cluster size and migration speed in a three-dimensional tissue environment, with sufficiently large clusters becoming mechanically hindered \cite{cai2016modeling}.

Our results are consistent with the idea that maintaining cluster size within a restricted developmental range may help preserve efficient migration through the confined nurse-cell environment. However, the present model does not establish that border-cell size or nurse-cell arrangement is optimized specifically for migration efficiency. Local nurse-cell geometry itself contributes to migration-path selection and varies along the migration route \cite{dai2020tissue,george2025chemotaxis}. We therefore interpret the observed architecture as a developmentally generated environment within which border-cell migration must remain robust, rather than as a geometry uniquely optimized for migration.

\subsection{Robustness to Geometry and Parameter Variations}

\begin{figure}[htp!]
    \centering
    \includegraphics[width=1\linewidth]{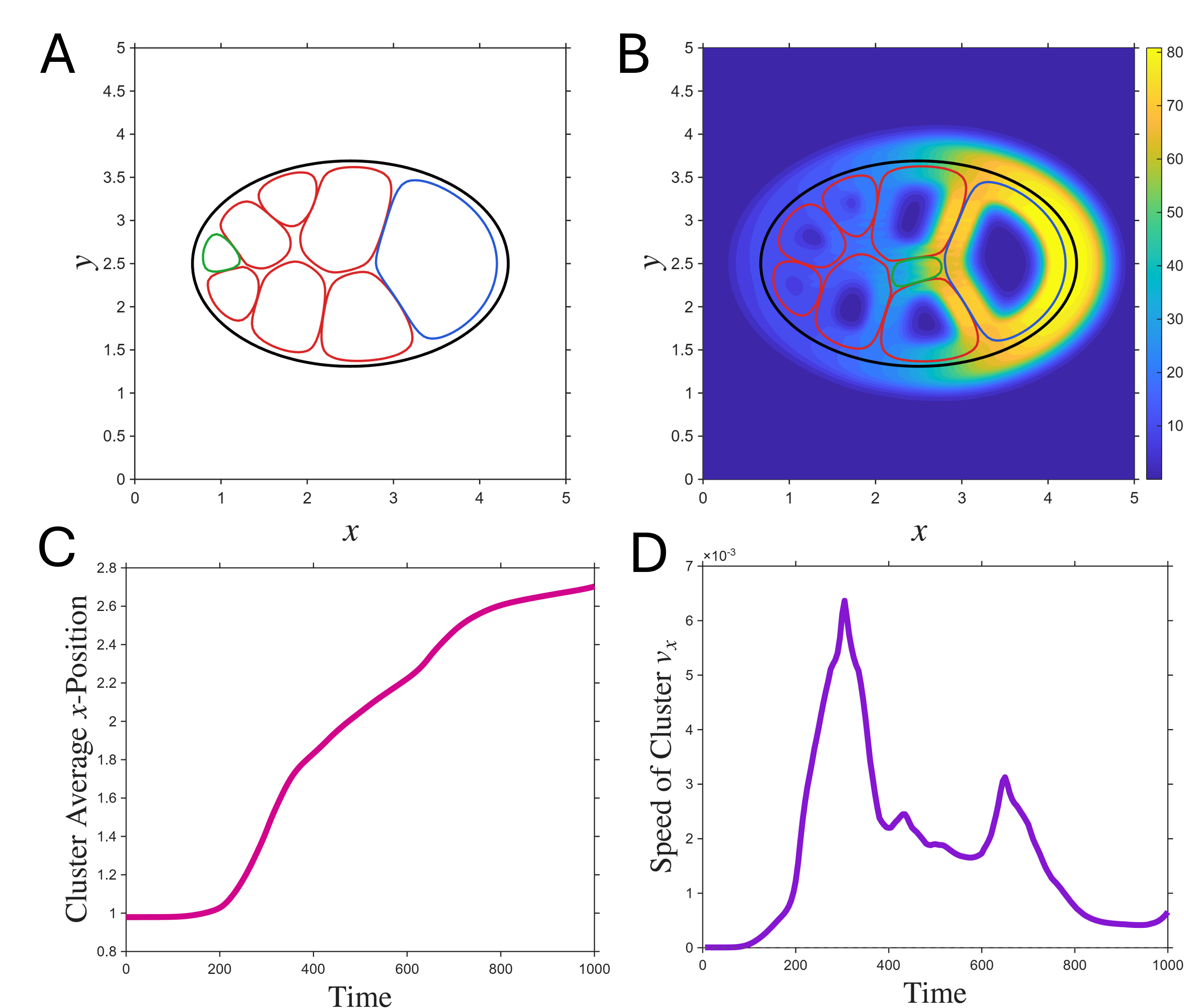}
    \caption{
     (A) An alternative
    egg-chamber geometry obtained by modifying the nurse-cell target volumes
    together with the initial nurse-cell placement. (B) corresponding
    concentration snapshot at $t=700$ with cell boundaries overlaid. (C) cluster average position as a function of time. (D)
    cluster migration speed as a function of time. In contrast to changing
    only the initial nurse-cell centers or initial nurse-cell sizes, this
    combined modification produces a genuinely different relaxed chamber
    while preserving a comparable overall oocyte--cluster arrangement. The
    resulting migration dynamics show two distinct slow-down regions,
    reflecting the two distinct gaps encountered by the cluster during
    migration.
    }
    \label{fig:comment9_geometry}
\end{figure}

{\allblack
Because the baseline simulations use a single representative nurse cell
arrangement and a fixed set of transport parameters, we examined whether the
main qualitative conclusions persist under changes in tissue geometry and
key parameters of the chemoattractant model.

% \begin{figure}[htp!]
%     \centering
%     \includegraphics[width=1\linewidth]{Figs/Figure10_main.png}
%     \caption{
%      (A) An alternative
%     egg-chamber geometry obtained by modifying the nurse-cell target volumes
%     together with the initial nurse-cell placement. (B) corresponding
%     concentration snapshot at $t=700$ with cell boundaries overlaid. (C) cluster average position as a function of time. (D)
%     cluster migration speed as a function of time. In contrast to changing
%     only the initial nurse-cell centers or initial nurse-cell sizes, this
%     combined modification produces a genuinely different relaxed chamber
%     while preserving a comparable overall oocyte--cluster arrangement. The
%     resulting migration dynamics show two distinct slow-down regions,
%     reflecting the two distinct gaps encountered by the cluster during
%     migration.
%     }
%     \label{fig:comment9_geometry}
% \end{figure}

We first considered an alternative egg chamber geometry. We modified the
 target areas of nurse-cells together with their initial placement to generate a
distinct chamber while preserving the overall
anterior--posterior arrangement of the border-cell cluster and oocyte.
Figure~\ref{fig:comment9_geometry} shows the resulting geometry,
chemoattractant distribution, cluster position, and migration speed. Directed
migration toward the oocyte is maintained, but the detailed speed profile
differs from that of the baseline geometry. In particular, two pronounced
slow-down regions occur as the cluster encounters two distinct confined gaps
along its migration path. Thus, the precise migration dynamics depend on
local nurse-cell packing, while the association between geometric
confinement and nonuniform migration speed persists in the alternative
configuration.

We next examined the sensitivity of the chemoattractant distribution to the
diffusivity transition parameter $s$ and secretion strength $\sigma$. Varying
$s$ changes the sharpness of the transition between low- and
high-diffusivity regions and consequently alters the contrast between
concentration levels inside and outside cellular domains. Nevertheless, the
overall posterior--anterior organization of the concentration field is
retained across the tested values (Supplementary Fig.~S7). In contrast,
changing $\sigma$ primarily modifies the magnitude of the concentration
field. After normalization by $\sigma$, the concentration profiles nearly
collapse, indicating that the secretion-strength parameter predominantly
rescales the concentration while preserving its overall spatial structure
(Supplementary Fig.~S8).

Finally, we tested the sensitivity of the boundary-localized oocyte source
to the diffuse-interface thickness by considering
$\epsilon_{\rm oct}^{2}=0.0005$, $0.001$, and $0.002$. The resulting
steady concentration fields retain the same overall spatial organization,
with only modest changes near the oocyte interface
(Supplementary Fig.~S9). Cross-sectional concentration profiles show
similarly small differences, while the normalized integrated source
strength, ${Q_{\epsilon}}/{Q_{\epsilon_0}}$, $Q_{\epsilon}=\sigma\int_{\Omega} \|\nabla\phi_{\rm oct}\|\,d\mathbf{x}$
% \[
%     \frac{Q_{\epsilon}}{Q_{\epsilon_0}},
%     \qquad
%     Q_{\epsilon}
%     =
%     \sigma\int_{\Omega}
%     \|\nabla\phi_{\rm oct}\|\,d\mathbf{x},
% \]
remains close to unity across the tested interface widths
(Supplementary Fig.~S10). This indicates that changing the resolved
diffuse-interface thickness modifies the local representation of the source
only modestly while preserving the total secretion and the large-scale
concentration field.
}

\section{Discussion}

Collective cell migration emerges from the interplay of molecular guidance cues and physical interactions with the tissue environment. In the \emph{Drosophila melanogaster} egg chamber, border-cell migration toward the oocyte is largely driven by chemotactic signaling, yet the spatial distribution of chemoattractant is also influenced by the crowded and heterogeneous geometry of the tissue. In this study, we build on our previously established phase-field representation of the egg chamber \cite{akhavan2026phase} to develop a spatially resolved model of chemoattractant {\allblack free diffusion} and couple this model with the Tangential Interface Migration (TIM) force to examine how chemical guidance, tissue geometry, and contact-mediated migration interact.

{\allblack
Alternative modeling approaches have also provided important insight into border-cell migration. In particular, Dai et al.~\cite{dai2020tissue} combined experimental observations with a coarse-grained physical model to demonstrate how tissue microtopography, adhesion, and chemotactic information contribute to migration-path selection. We view such biophysical and phase-field descriptions as complementary rather than mutually exclusive.
The phase-field representation is particularly useful for the questions considered here because cell boundaries are represented as continuous spatial fields and can therefore enter directly into both the chemical and mechanical components of the model. In the present framework, the cellular phase fields determine the geometry-dependent diffusion coefficient,
localize secretion to the oocyte interface, and identify cluster--nurse-cell contact regions and tangential directions used by the
TIM force. Thus, the same interface representation couples tissue geometry, extracellular distribution, and contact-mediated migration without requiring explicit boundary tracking. We do not claim that the qualitative behaviors described here are intrinsically irreproducible by traditional approaches; rather, the phase-field formulation provides a convenient framework for resolving and coupling spatially distributed chemical fields with deformable cellular interfaces.
}

A key contribution of this work is the formulation of a diffusion--reaction model in which chemoattractant is secreted along the oocyte boundary and diffusion is modulated by the cellular geometry.
{\allblack
The phase-dependent diffusion coefficient strongly reduces free movement within cellular regions while permitting greater diffusion through the extracellular spaces between cells.}
This view is consistent with broader work on morphogen transport emphasizing that signaling molecules move through geometrically complex and tortuous tissue environments, which can substantially influence gradient formation and interpretation \cite{muller2013morphogen}.
{\allblack
The simulations consequently produce heterogeneous and anisotropic concentration distributions whose local structure depends on the arrangement of nurse cells and the available extracellular corridors. The association between migration-speed reductions and geometrically confined regions is retained in an alternative nurse-cell packing, although the precise speed history depends on the particular geometry. The parameter sensitivity
analyses likewise show that quantitative concentration contrast depends on the diffusivity-transition parameter $s$, whereas changing $\sigma$ primarily rescales concentration magnitude and varying the resolved oocyte interface thickness produces comparatively modest changes in the integrated source and large-scale concentration field. These results support the robustness of the qualitative geometry-dependent chemoattractant distribution while also demonstrating that its quantitative expression depends on model parameters and tissue architecture.
}
Our findings are consistent with experimental observations showing that local tissue geometry and topography influence border-cell migration and path selection, as well as with our previous observations that tissue architecture modulates chemoattractant dispersion along the migration route
\cite{dai2020tissue,george2025chemotaxis}.

{\allblack
The dynamic--static comparisons further clarify the specific role of the time-dependent chemoattractant formulation. When the chemoattractant field is initialized from $c_0=0$, its transient establishment produces a substantial delay in the onset and progression of migration relative to simulations initialized from an already established steady-state field. By contrast, when both simulations begin from the same $c_{\mathrm{ss}}$, continuously updating the concentration field and holding that field fixed produce much more similar migration trajectories, with smaller quantitative differences at later times. {\allblack However, it is notable that dynamic change in chemaottractant to the landscape tended to increase the speed with which the cluster was able to cross the chamber.} {\allblack Thus, in the parameter regime examined here, the strongest effect of the time-dependent diffusion--reaction model occurs during establishment of the guidance field; a static geometry-shaped steady-state profile provides a reasonable approximation for much of the subsequent migration once that field has already formed.} }

{\allblack Nevertheless, the dynamic formulation captures a feedback that a fixed concentration profile cannot represent. The local probe analysis shows that motion of the cluster and rearrangement of neighboring nurse-cell interfaces are accompanied by changes in local chemoattractant concentration when the field is allowed to evolve, whereas the corresponding static concentration remains unchanged by construction. The evolving cellular geometry therefore modifies both the mechanical confinement experienced by the cluster and the phase-dependent diffusion landscape through which the chemoattractant is transported. Accordingly, local migration-speed variations should not be attributed exclusively to either chemical or mechanical effects; the model instead provides a framework in which both consequences of tissue geometry can act simultaneously.
}

To examine how these chemical cues contribute to migration, we couple the concentration field to the TIM force \cite{akhavan2026phase}, which represents contact-mediated tangential traction at interfaces between the border-cell cluster and surrounding nurse cells.
{\allblack
In the present formulation, the local receptor-response function controls the magnitude of chemically activated TIM traction, while the projected concentration gradient determines its tangential direction. The use of a sign function therefore represents an idealized directional response and does not explicitly include a gradual loss of directional sensitivity in very shallow gradients. Likewise, the deterministic formulation does not include receptor noise or a minimum detectable gradient and therefore does not predict the fundamental size-dependent limits of chemotactic sensing. These features provide natural directions for extending the model when quantitative information on sensing precision becomes available.
}

{\allblack
Cluster-size simulations further demonstrate that the size of the migrating group affects both components of the coupled model. Changing the target cluster area modifies the timing and magnitude of migration and also changes the evolving local chemoattractant concentration through the geometry-dependent diffusivity. Border-cell number is developmentally regulated through the signaling network that specifies migratory follicle
cell fate \cite{monahan2013socs36e,montell2012group,hombria2002fertile},
and previous work has shown that sufficiently large clusters can become mechanically hindered \cite{cai2016modeling}. Our results are therefore consistent with the idea that maintaining cluster size within a restricted developmental range can support robust collective migration. However, they do not establish that either cluster size or nurse cell arrangement is uniquely optimized for migration efficiency. Rather, the cluster migrates within a developmentally generated tissue architecture in which collective organization may help preserve directional information while excessive size increases confinement and mechanical resistance.
}

{\allblack
The present simulations should also be interpreted as mechanistic rather than quantitatively fitted reconstructions of the \emph{in vivo} chemoattractant field. The PVF1-eGFP imaging experiments motivate localization of the source and provide qualitative information about ligand distribution, but they are not used to calibrate concentration magnitude or infer model parameters. Similarly, the target cellular areas and the effective adhesion, repulsion, and interfacial parameters are not independently inferred from image segmentation or mechanical measurements. The latter are coarse-grained model parameters selected within the regime of the previous phase-field formulation to maintain stable and representative cellular morphologies and interactions. Quantitative image-based estimation of cell geometry and calibration or inference of effective mechanical and transport parameters would therefore be an important extension of the present framework, particularly as suitable spatially resolved experimental data become available.
}

{\allblack
As an additional numerical check, comparison with a sharp-interface finite-element formulation produces the same overall organization of the chemoattractant field and qualitatively similar posterior migration of the cluster, although quantitative differences remain in concentration, instantaneous speed, and final position. This supports the conclusion that the principal geometry-dependent behavior is not solely an artifact of the diffusion-interface representation, while also illustrating the computational advantage of the phase-field formulation for simulations involving evolving cellular boundaries.
}

{\allblack
An important limitation of the present framework is its two-dimensional representation of the egg chamber. In three dimensions, extracellular ligand can access additional pathways around neighboring cells, and therefore the magnitude and local shape of geometry-induced concentration gradients may differ from those obtained in a planar cross-section. Likewise, cell--cell contacts become surfaces rather than curves, providing additional directions for cluster deformation and tangential motion. We therefore interpret the present results as demonstrating qualitative mechanisms by which cellular geometry can redirect extracellular signal distribution and modulate contact mediated migration, rather than as quantitative predictions of three-dimensional chemoattractant concentrations, mechanical stresses, or migration speeds. Extending the coupled chemo-mechanical framework to experimentally reconstructed three-dimensional geometries will be an important direction for future work.
}{\allblack As we refine our models, it will also be worthwhile to incorporate uptake of signal by the migrating cells themselves, to examine how this may modulate local chemoattractant gradients in limiting conditions.}

In summary, our results show that tissue geometry influences border-cell migration through both physical confinement and its effect on chemoattractant transport.
{\allblack
The time-dependent chemoattractant field is particularly important during transient gradient establishment, while moving interfaces continue to produce measurable local redistribution of the chemical field during migration.}
The integrated framework therefore provides a means of studying how deformable tissue architecture, extracellular signaling, and
contact-mediated forces interact during collective cell migration.

\section{Methods}

\subsection{Phase Field Modeling of the Egg Chamber}

To model the geometry and cellular architecture of the \emph{Drosophila melanogaster} egg chamber, we employ a multiphase-field framework \cite{nonomura2012study, lee2017new, seirin2021extra, moure2021phase} in which each biological region (namely, the oocyte, border cell cluster, and individual nurse cells) is represented by a scalar-valued phase field variable \( \phi_m(x, y, t) \in [0, 1] \). In this representation, \( \phi_m \approx 1 \) corresponds to the interior of the \( m \)-th domain and \( \phi_m \approx 0 \) to the exterior. For the epithelial layer ($\phi_0$), we take the phase field to be 1 inside (and beyond - external to the chamber) and 0 inside the chamber.
{\allblack The present model represents a two-dimensional cross-section of the egg chamber rather than a full three-dimensional reconstruction.
A biological egg chamber contains 15 nurse cells and one oocyte, whereas the six nurse-cell domains included here constitute a reduced two-dimensional representation of the local cellular architecture surrounding the border-cell migration path rather than the complete nurse-cell population \cite{montell2012group, peercy2020clustered}.
The phase-field formulation itself is not restricted to two dimensions and has been applied to multicellular systems in both two and three dimensions \cite{nonomura2012study, akiyama2019numerical, moure2021phase}. Here, the two-dimensional representation allows us to resolve cell interfaces, extracellular corridors, and contacts between the border cell cluster and neighboring nurse cells while maintaining a tractable framework for coupling cell mechanics with spatial chemoattractant dynamics. Accordingly, the model is intended primarily to examine qualitative mechanisms arising from geometry-dependent transport and cell--cell interactions rather than to provide a complete quantitative representation of three-dimensional egg-chamber mechanics.}
% \sout{however, we adopt the opposite convention: the phase field is taken to be approximately \(0\) inside the epithelial region and approximately \(1\) outside it.} }
This diffuse interface approach avoids the need for explicit boundary tracking and naturally accommodates dynamic shape changes, migration, and domain interactions through smooth transitions across interfaces.
The time evolution of each phase variable \( \phi_m \) is governed by \cite{nonomura2012study, akhavan2026phase}:
\begin{align}
    \frac{\partial \phi_m}{\partial t} = -\mu \frac{\delta E}{\delta \phi_m}, \quad (1 \le m \le N), \label{wrt phi_m} 
\end{align}
where \( \mu \) is a mobility parameter, \( N \) is the total number of phase field variables, and \( E \) is the total free energy of the system.
The specific form of the governing equation is: 
\begin{align}\label{eq: time evol}
    -\frac{1}{\mu} \frac{\partial \phi_m}{\partial t} = \epsilon_m^2 \nabla^2 \phi_m + \phi_m(1-\phi_m) \left[\phi_m - \frac{1}{2} + 6 \mathcal{F}(\phi_m, \phi_0)\right],
\end{align}
where the first term promotes smooth interface curvature and the second term defines a double-well potential with coupling to a system-wide functional \( \mathcal{F} \).
The functional \( \mathcal{F}(\phi_m, \phi_0) \) encodes mechanical and spatial constraints:
\begin{align}\label{f function}
\begin{split}
    \mathcal F(\phi_m, \phi_0) =& 2  \alpha_m(V_m(t) - \bar V_m(t)) + \beta_0 h(\phi_0) +\beta(m,n) (\xi -h(\phi_m)) \\
    &+2 \alpha_0 \left[\int (1-h(\phi_0))dx - \sum_{m=1}^N V_m)\right]+\gamma_0 \nabla^2h(\phi_0)+\gamma(m,n) \nabla^2(\xi -h(\phi_m)).
    \end{split}
\end{align}
% where  $\phi_0$ is the phase variable of epithelial layer and \( h(\phi) = \phi^2 (3 - 2\phi) \) is a smooth interpolation function, \( V_m(t) = \int h(\phi_m) \, dx \) is the volume of the \( m \)-th region, \( \bar{V}_m(t) \) is its target volume, and \( \xi = \sum_{m=1}^N h(\phi_m) \) tracks the overall occupancy of the domain.
{\allblack Here, $(\phi_0)$ is the prescribed phase field representing the epithelial boundary, with $(\phi_0=0)$ in the interior of the egg chamber and $(\phi_0=1)$ in the epithelial and exterior region. The epithelial phase field is held fixed during the simulations and therefore acts as the confining boundary within which the cellular phase fields evolve. The interpolation function $(h(\phi)=\phi^2(3-2\phi))$ provides a smooth approximation of the cellular indicator function, $(V_m(t)=\int_\Omega h(\phi_m)d\mathbf{x})$ is the volume (area in the present two-dimensional model) of the $(m)$-th cellular domain, $(\bar V_m(t))$ is its target volume, and $(\xi=\sum_{m=1}^{N}h(\phi_m))$ measures the total occupancy of the cellular phase fields.

The terms in $(\mathcal F)$ represent distinct mechanical constraints. The term $(2\alpha_m(V_m-\bar V_m))$ arises from the individual volume penalty for cellular domain (m) and maintains its volume near the prescribed target. The coefficient $(\beta_0)$ controls exclusion between cellular phases and the prescribed epithelial region, whereas $(\beta(m,n))$ controls exclusion between neighboring cellular domains and prevents nonphysical cell overlap. The term proportional to $(\alpha_0)$ does not control the volume of the epithelial layer; rather, it arises from a global interior-volume constraint that penalizes mismatch between the area enclosed by the epithelial boundary, $(\int_\Omega(1-h(\phi_0)),d\mathbf{x}$), and the total area occupied by the modeled cellular domains. Finally, $(\gamma_0)$ and $(\gamma(m,n))$ describe cell--epithelial and cell--cell adhesive interactions, respectively. The diffuse-interface width and regularization of each cellular phase are controlled separately by $(\epsilon_m)$ in Eq.~\eqref{eq: time evol}. A full derivation of these energy contributions is provided in our previous phase-field formulation \cite{akhavan2026phase}.
}

The energy terms enforce volume conservation (\( \alpha_m \)), exclusion between overlapping phase fields (\( \beta(m,n) \)), and interface regularity (\( \gamma(m,n) \)). The shared extracellular domain is implicitly defined as the region where all \( \phi_m \ll 1 \), allowing us to treat extracellular space and intracellular domains within the same computational framework.

This modeling structure, originally developed in our prior study~\cite{akhavan2026phase}, enables robust simulation of cell migration within crowded tissue. In the present work, we build on this framework to examine how extracellular geometry, defined by the collective layout of the phase fields, shapes the distribution of chemoattractant concentration throughout the egg chamber.
We represent the oocyte, border cell cluster, and surrounding nurse cells using separate phase field variables, which together define the equilibrium architecture of the egg chamber (Fig.~\ref{fig:equilibrium}(A)).

% We use a phase-field approach to represent the geometry of the \emph{Drosophila} egg chamber, including the oocyte, nurse cells, and the migrating border cell cluster. This continuum-based framework allows us to simulate interfaces between cellular domains and extracellular space without explicitly tracking their boundaries. Each domain is described by a smooth phase field variable \( \phi_i(x,y,t) \in [0,1] \), where \( \phi_i \approx 1 \) denotes the interior of the \( i \)-th cell or region, and \( \phi_i \approx 0 \) corresponds to the surrounding medium. 

% The dynamics of each \( \phi_i \) are governed by Allen–Cahn-type evolution equations with interface-penalizing terms to ensure stable, compact domains. We impose volume conservation constraints to reflect biological incompressibility of cellular compartments. The nurse cells are arranged in a compact lattice, while the oocyte is modeled as a large posterior domain. The border cell cluster is initialized as a small, mobile domain at the anterior side, able to migrate through the nurse cell region.

% \newpage
\subsection{Spatial Modeling of Chemoattractant Concentration}
In developing tissues, chemoattractant gradients play a pivotal role in guiding the collective migration of cells. One well-characterized example is the directional movement of the border cell cluster toward the oocyte in response to spatially regulated chemoattractant cues in the \textit{Drosophila} egg chamber \cite{duchek2001guidance, mcdonald2006multiple, prasad2007cellular, montell2012group}. To elucidate the mechanisms underlying gradient formation in such a constrained and heterogeneous environment, we developed a diffusion-reaction model in which the effective diffusion of the chemoattractant is modulated by the surrounding tissue architecture.

\subsubsection{One-Dimensional Illustrative Model of Phase-Dependent Diffusion}
In the \emph{Drosophila} egg chamber, chemoattractant ligands, including PVF1, are produced at the oocyte and diffuse through narrow extracellular corridors to guide the migration of the border cell cluster. 
% {\allcblue (probably no evidence on it! - \textbf{there is some evidence} - cite Duchek for pvf1, Dai/topography paper for Krn, our paper) }
Importantly, these signaling molecules seem to be secreted by the oocyte~\cite{duchek2001guidance, duchek2001guidance2} and are not appreciably found in the interiors of nurse cells~\cite{ dai2020tissue, george2025chemotaxis}. 
% {\allcpink (could also note here, that even when overexpressed, appear restricted to oocyte as in figure 3) (also, just a note to you- it is found in border cells but only when bound to receptor- so i took that out here)}
Instead, they are largely confined to extracellular regions, where their distribution is shaped by the geometry of the surrounding tissue \cite{mekus2018effects, george2025chemotaxis}. Experimental imaging suggests that chemoattractant gradients are influenced by non-responsive cells in the domain and may be distorted or weakened in regions where extracellular space is restricted \cite{dai2020tissue, george2025chemotaxis}, 
% {\allcpink (i think we should cite Dai here as it relates to krn but may had to alter the statement- see if you like it)}
suggesting that spatial heterogeneity in tissue structure plays a critical role in modulating the effective guidance cues available to migrating cells.

{\allblack Before introducing the full two-dimensional chemoattractant model of the egg chamber, we first consider a reduced one-dimensional problem to isolate the effect of a moving low-diffusivity phase-field domain on chemoattractant transport. This one-dimensional model is not the concentration model used in the coupled egg-chamber simulations and is not obtained from a fixed-$y$ cross-section of the two-dimensional system. Rather, it serves as a simplified mechanistic illustration of how phase-dependent diffusivity can exclude chemoattractant from a cellular domain and reshape the local concentration profile as that domain moves.}

% To examine how the geometry and motion of a low-diffusivity domain affects the accumulation and profile of chemoattractant over time, we construct a one-dimensional model in which diffusion is restricted by a moving, phase-field-defined region that mimics a migrating cluster.

% {\allblack For this purpose, we construct a one-dimensional model in which diffusion is modulated by a moving phase-field-defined region that represents a low-diffusivity cellular domain.}

\noindent {\bf{Diffusion–Reaction Model with Phase-Dependent Diffusivity.}}  
We model the chemoattractant concentration \( c(x,t) \) using a diffusion–reaction equation with diffusion restricted to the domain outside of phases:
\begin{align}
    \frac{\partial c}{\partial t} = \frac{\partial }{\partial x} \left(D(\phi(x)) \frac{\partial c}{\partial x} \right) - kc, \label{eq: diff pde}
\end{align}
where \( c(x,t) \) denotes the concentration of the chemoattractant at position \( x \) and time \( t \); \( D(\phi) \) is a phase-dependent diffusion coefficient; \( k \) is the decay rate; and \( \phi(x) \) is a smooth phase field variable that represents the presence of intracellular regions. 

To restrict diffusion within cell interiors while allowing transport in extracellular space, we define the diffusion coefficient as:
\begin{align}\label{eq: diffusion}
    D(\phi) = \frac{D_0}{1 + e^{(\phi - \phi_*)/s}},
\end{align}
where \( D_0 \) is the maximal diffusion constant, \( \phi_* \) is a threshold indicating the boundary between intra- and extracellular domains, and \( s \) controls the sharpness of this transition. 
{\allblack This sigmoidal form provides a smooth transition between high- and low-diffusivity regions while retaining a diffuse-interface representation of the cellular domain~\cite{stark2023open}.}
% This sigmoidal form enables spatially graded but biologically realistic restriction of diffusivity, as commonly used in phase field models of heterogeneous media

\begin{figure}[!hpt]
    \centering
    \includegraphics[width=1\linewidth]{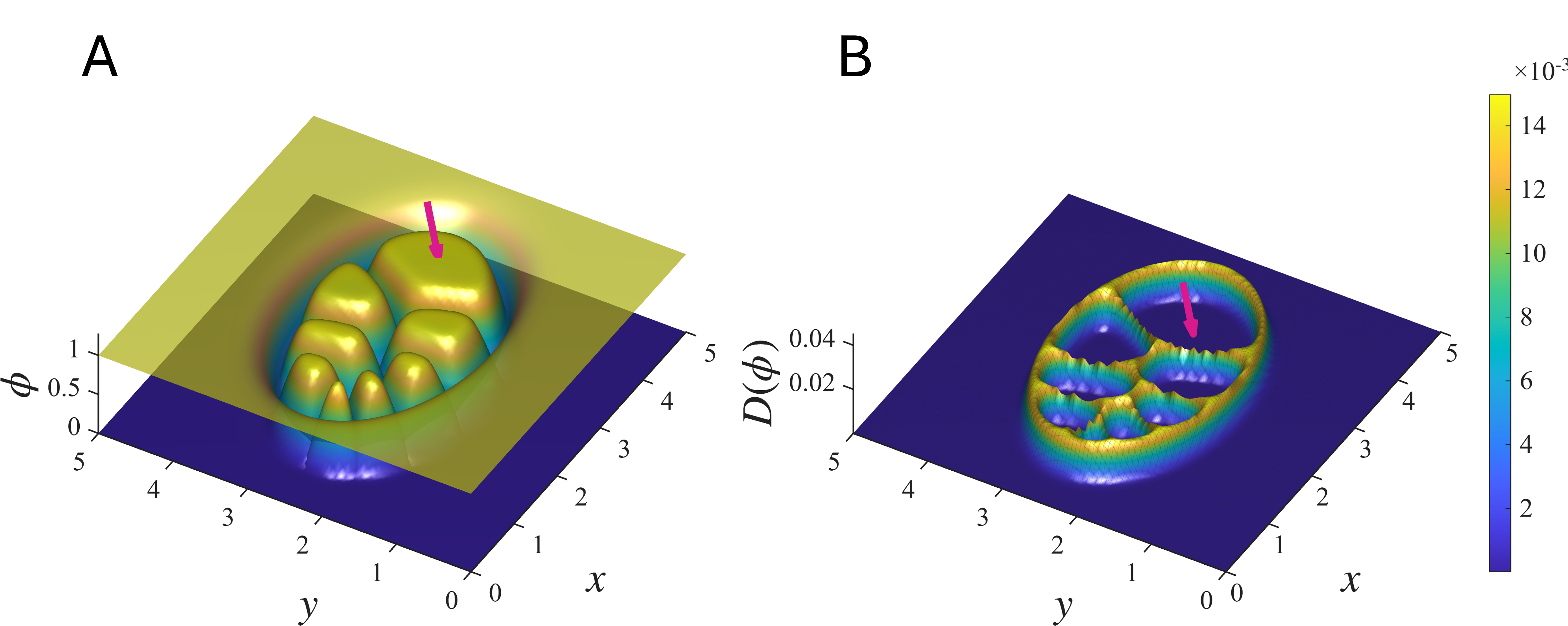}
\caption{{\bf{Architecture of \emph{Drosophila} egg chamber.}} 
% (A) Phase-field representation of the egg chamber structure, showing six nurse cells, the border cell cluster positioned at the anterior, and the oocyte at the posterior. 
{\allblack (A) Reduced two-dimensional representation of the egg-chamber architecture used in the model. The border cell cluster is initially located at the anterior and migrates through the nurse-cell environment toward the posterior oocyte. Six nurse-cell phase-field domains are included to represent the local cellular geometry surrounding the migration path; these do not represent the full biological complement of 15 nurse cells. The surrounding epithelial phase field defines the confining boundary of the modeled egg chamber, while regions between the cellular domains represent extracellular space.}
(B) Spatial profile of the diffusion coefficient \(D(\phi)\) (Eq.~\eqref{eq pfm diff}), {\allblack illustrating reduced diffusivity within cellular regions and higher diffusivity in the extracellular space. The magenta arrow indicates the oocyte in both panels. Parameter values are provided} in Supplementary Table~1.}
    \label{fig:equilibrium}
\end{figure}

The corresponding diffusion profile of Fig.~\ref{fig:equilibrium}(A) is shown in Fig.~\ref{fig:equilibrium}(B), illustrating the spatially varying  effective diffusion coefficient ($D(\phi)$) (See also Supplementary video S4 and S6).

{\allblack 
\noindent {\bf Prescribed moving phase-field profile.}
To represent a moving low-diffusivity cellular domain in the one-dimensional illustrative model, we prescribe a smooth phase-field profile constructed from two opposing hyperbolic-tangent interfaces. Such $\tanh$ profiles are commonly associated with diffuse-interface and Allen--Cahn formulations \cite{allen1979microscopic,moure2021phase}. Here, however, the Allen--Cahn dynamics are not evolved; the profile is used only as a convenient smooth representation of a translating cellular domain.

Specifically, we define
\begin{align}\label{eq:tanh_window}
    \phi(x)=\tfrac12\!\left[
    \tanh\!\Big(\frac{x-(x_0-R)}{\sqrt{2D_\phi}\,\varepsilon}\Big)
    -\tanh\!\Big(\frac{x-(x_0+R)}{\sqrt{2D_\phi}\,\varepsilon}\Big)
    \right],
\end{align}
where $x_0$ denotes the center of the phase-field domain, $R$ is its half-width, and $\varepsilon$ controls the diffuse-interface thickness. The parameter $D_\phi$ sets the spatial scaling of the interface. Thus, $\phi\approx1$ on $[x_0-R,x_0+R]$ and smoothly approaches $0$ outside this interval. To prescribe motion with constant speed $v$, the center is translated according to $x_0(t)=x_0(0)+vt$, equivalently yielding a profile of the form $\phi(x-vt)$.

This construction is used only in the reduced one-dimensional example to illustrate how motion of a low-diffusivity phase-field domain modifies the surrounding chemoattractant distribution.}
% To model motion at speed $v$, we translate the mask kinematically as $\phi(x-vt)$.
% {\allcblue (Maybe a bit too much explanation here??? MSG thinks it is okay)}{\allorange BEP agrees that it is okay.  Like it.}

\begin{figure}[!htp]
    \centering
    \includegraphics[width=1\linewidth]{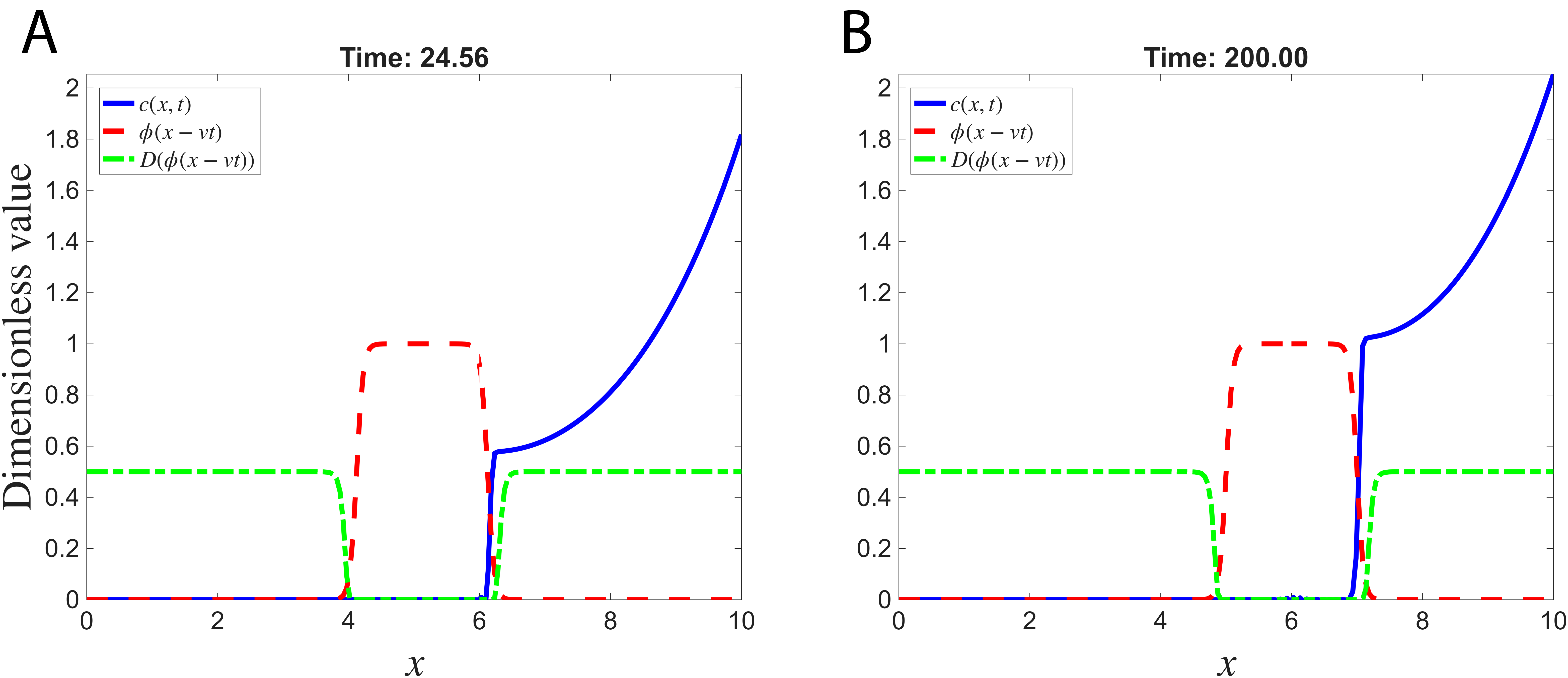}
    \caption{{\bf One-dimensional illustration of phase-dependent chemoattractant transport. }Evolution of the chemoattractant concentration \(c(x,t)\) (blue curve), the {\allblack moving phase-field profile} \(\phi(x-vt)\)  (red, dashed curve), and the {\allblack corresponding} diffusion coefficient \(D(\phi(x-vt))\) (green, dot-dashed curve) over time. 
    Panel (A) shows the system at an early stage, while panel (B) shows the translated phase-field domain at the end of the simulation. {\allblack The horizontal axis denotes position $(x)$}. The parameters used are $(R=1)$ and $(\varepsilon=0.07)$; all remaining parameter values for the one-dimensional illustrative model are listed in Supplementary Table~1 (see also in  Supplementary video S1).}
    \label{fig:Allen Cahn solutions}
\end{figure}

Figure~\ref{fig:Allen Cahn solutions} shows the evolution of the system over time. Panel (A) shows the initial condition where the phase field variable $\phi(x)$ is centered around $x=5$, while panel (B) shows its position at $x=6$ after migrating with constant velocity. 
The diffusion coefficient $D(\phi(x))$ sharply decreases within the region defined by $\phi \approx 1$, restricting diffusion to the exterior. As a result, the chemoattractant concentration $c(x,t)$ accumulates near the leading edge of the moving region and remains excluded from the interior. 

% This behavior demonstrates that coupling phase-field geometry to a diffusion coefficient effectively confines the chemoattractant to extracellular domains 
% % {\allblack [I just realized that the total chemoattractant appears to be conserved, which I hadn't noticed before - not sure it needs a comment but a reviewer may note that]}
% , reproducing a key biophysical constraint observed \emph{in vivo}.
% Our model provides a continuum-level mechanism for studying signal confinement during migration and offers a computational framework for examining how geometric barriers contribute to the formation of asymmetric, biologically relevant chemical gradients.
{\allblack
This behavior illustrates how coupling a low-diffusivity phase-field domain to the diffusion coefficient can restrict chemoattractant transport through that domain and modify the surrounding concentration profile. The purpose of this reduced example is to isolate this transport mechanism in a simple setting before considering the geometrically complex two-dimensional egg-chamber model.
}

% \begin{figure}[H]
%     \centering
%     \includegraphics[width=1\linewidth]{Figs/Figure2_main.png}
%     \caption{{\bf One-dimensional illustration of phase-dependent chemoattractant transport. }Evolution of the chemoattractant concentration \(c(x,t)\) (blue curve), the {\allblack moving phase-field profile} \(\phi(x-vt)\)  (red, dashed curve), and the {\allblack corresponding} diffusion coefficient \(D(\phi(x-vt))\) (green, dot-dashed curve) over time. 
%     Panel (A) shows the system at an early stage, while panel (B) shows the translated phase-field domain at the end of the simulation. {\allblack The horizontal axis denotes position $(x)$}. The parameters used are $(R=1)$ and $(\varepsilon=0.07)$; all remaining parameter values for the one-dimensional illustrative model are listed in Supplementary Table~1 (see also in  Supplementary video S1).}
%     \label{fig:Allen Cahn solutions}
% \end{figure}
{\allblack
The reduced one-dimensional example therefore illustrates how motion of a low-diffusivity domain can locally redistribute chemoattractant concentration. In the following section, we formulate the full two-dimensional chemoattractant model independently, using the phase-field geometry of the egg chamber to determine spatially varying diffusion and oocyte-localized secretion.
}
% The combination of a traveling tanh-profile phase field and nonlinear diffusion leads to persistent chemical accumulation ahead of the migrating domain. This profile sets up a spatial gradient that can, in turn, be used to guide the motion of cells in higher-dimensional models, as explored in later sections.

\subsubsection{Chemoattractant Concentration in the Egg Chamber}

\begin{figure}[!hpt]
    \centering
    \includegraphics[width=1\textwidth]{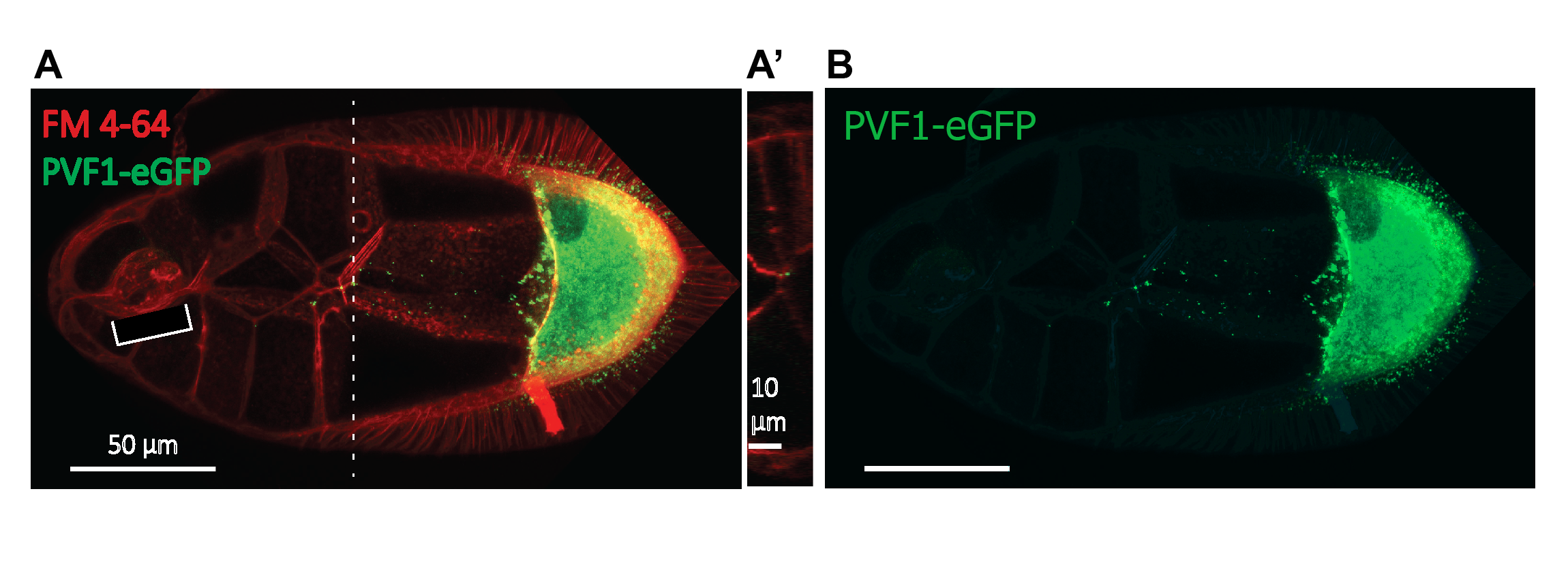}
    \caption{Live Distribution of Overexpressed PVF1-eGFP. (A) Representative egg chamber expressing PVF1-eGFP (green) under the control of mat$\alpha$-Gal4.
     Egg chamber was labeled with lipophilic FM4-64 dye (red) to label cell membranes and imaged live in a confocal microscope; maximum intensity projection is shown in A and the dotted line is shown in cross section in (A').  Bracket denotes the border cell cluster.  (B) shows PVF-eGFP alone; note the that the majority of the signal is detected in the oocyte, but some is secreted and detectable along the border cell migration path. 
     % {\allred (Caption is updated)}
     }
    \label{fig:eggchamberimages}
\end{figure}

Live imaging {\allblack experiments performed using the} 
% of developing \emph{Drosophila} egg chambers in which a
tagged chemoattractant, PVF1-eGFP \cite{george2025chemotaxis}
% {\allcpink cite George here for making this line}
over-expressed in the germline cells
{\allblack show strong enrichment of PVF1-eGFP in the oocyte, together with detectable signal along the oocyte surface and within the border-cell migration region (Fig.~\ref{fig:eggchamberimages}). 
The images shown in Fig.~\ref{fig:eggchamberimages} are representative images from these experiments and are distinct from those presented in our previous study \cite{george2025chemotaxis}, in which the spatial distribution of PVF1-eGFP and its relationship to tissue geometry were characterized experimentally. These observations are consistent with previous evidence that PVF1 is produced by the oocyte and provides a guidance cue for border-cell migration \cite{duchek2001guidance,mcdonald2003pvf1}.
The fluorescence images are used here to motivate the spatial localization of the chemoattractant source and are not used for quantitative fitting or calibration of the simulated concentration field.}

{\allblack Motivated by these observations, we make the simplifying modeling assumption that chemoattractant secretion is distributed along the oocyte boundary. We do not interpret the fluorescence pattern in Fig.~\ref{fig:eggchamberimages} as demonstrating a spatially uniform  secretion rate; rather, the distributed boundary source provides a coarse-grained representation of ligand release from the oocyte surface. To incorporate this assumption into the phase-field framework, we represent the oocyte by the phase-field variable $\phi_{\mathrm{oct}}$ and define the chemoattractant concentration $c(x,y,t)$ to evolve according to}
\begin{align}\label{eq: chemo for oocyte}
    \frac{\partial c}{\partial t} = \nabla \cdot \left(D(\phi) \nabla c \right) - kc + \sigma \|\nabla \phi_{\text{oct}}\|,
\end{align}
% {\allblack where $D(\phi)$ is the spatially varying diffusion coefficient, $k$ is the chemoattractant degradation rate, and $\sigma$ controls the strength of chemoattractant secretion. Because $\|\nabla\phi_{\mathrm{oct}}\|$ is appreciable only across the diffuse interface of the oocyte phase field, the source term $\sigma\|\nabla\phi_{\mathrm{oct}}\|$ localizes chemoattractant production to the oocyte boundary. 

% As the phase-field interface becomes narrower, $\|\nabla\phi_{\mathrm{oct}}\|$ increases in magnitude while its spatial support decreases, such that the integrated source remains finite. We assume a constant value of $\sigma$ along the oocyte interface, corresponding to spatially uniform secretion per unit boundary length in the present two-dimensional model. This is a simplifying modeling assumption and is not intended to imply that the microscopic secretion rate is experimentally known to be uniform along the entire oocyte surface.
% {\allred
% The robustness of the diffuse-interface source with respect to the chosen interface width was verified numerically by varying $\varepsilon_{\rm oct}$ around its baseline value; the integrated source strength and resulting steady-state concentration profiles showed only minor variation (Supplementary Fig.~...).}
{\allblack
where $D(\phi)$ is the spatially varying diffusion coefficient, $k$ is the chemoattractant degradation rate, and $\sigma$ controls the strength of chemoattractant secretion. Because
$\|\nabla\phi_{\mathrm{oct}}\|$ is appreciable only across the diffuse interface of the oocyte phase field, the source term $\sigma\|\nabla\phi_{\mathrm{oct}}\|$ localizes chemoattractant production to the oocyte boundary.

Although the pointwise magnitude of $\|\nabla\phi_{\mathrm{oct}}\|$ depends on the phase-field interface thickness, its spatial support changes correspondingly. As the interface becomes narrower, $\|\nabla\phi_{\mathrm{oct}}\|$ increases in magnitude while the width of the source region decreases, so that the integrated source
\[
\sigma\int_{\Omega}\|\nabla\phi_{\mathrm{oct}}\|\,d\mathbf{x}
\]
remains finite and, for a sufficiently resolved diffuse interface,
approximates secretion proportional to the oocyte boundary length.

The coefficient $\sigma$ is taken to be spatially constant along the oocyte
interface. Thus, the present two-dimensional model assumes a spatially
uniform secretion rate per unit boundary length. This is a simplifying
modeling assumption: Fig.~\ref{fig:eggchamberimages} is used to motivate
localization of the source to the oocyte boundary, but the fluorescence
intensities are not used to infer either uniform  secretion or
the quantitative value of $\sigma$. Spatially heterogeneous secretion could
be incorporated in future extensions by allowing $\sigma$ to vary with
position along the oocyte boundary.

% {\allred 
% The robustness of the diffuse-interface source with respect to the oocyte interface thickness was also examined numerically by varying
% $\varepsilon_{\rm oct}$ around its baseline value. The resulting integrated source strength and chemoattractant concentration field showed only minor variation over the numerically resolved interface widths considered (Supplementary Fig.~).
% % The integrated source varied by less than [X]\%, while the relative
% % $L^2$ difference in the concentration field remained below [Y]\%.
% }

% For the present simulations, we impose homogeneous Neumann boundary conditions, \[ \frac{\partial c}{\partial n}=0, \] on the outer boundary of the computational domain. This condition is used as a simplifying approximation to isolate chemoattractant transport within the modeled egg-chamber cross-section over the timescale considered; it is not intended to imply that the biological follicular boundary is perfectly impermeable to molecular transport.
}

% where \( D(\phi) \) is the spatially varying diffusion coefficient, \( k \) is the decay rate of the chemoattractant, and \( \sigma \|\nabla \phi_{\text{oct}}\| \) represents a boundary-localized source term. The gradient norm \( \|\nabla \phi_{\text{oct}}\| \) is nonzero only near the interface of the oocyte’s phase field, thereby ensuring that secretion is confined to the oocyte surface.
% We also consider zero-flux (homogeneous Neumann) boundary conditions on all boundaries of the computational domain $\frac{\partial c}{\partial {\bf{n}}} = 0$.
% \begin{figure}[!hpt]
%     \centering
%     \includegraphics[width=1\textwidth]{Figs/Figure3_main.png}
%     \caption{Live Distribution of Overexpressed PVF1-eGFP. (A) Representative egg chamber expressing PVF1-eGFP (green) under the control of mat$\alpha$-Gal4.
%      Egg chamber was labeled with lipophilic FM4-64 dye (red) to label cell membranes and imaged live in a confocal microscope; maximum intensity projection is shown in A and the dotted line is shown in cross section in (A').  Bracket denotes the border cell cluster.  (B) shows PVF-eGFP alone; note the that the majority of the signal is detected in the oocyte, but some is secreted and detectable along the border cell migration path. 
%      % {\allred (Caption is updated)}
%      }
%     \label{fig:eggchamberimages}
% \end{figure}

The diffusion coefficient is also modified to account for multiple cellular domains—nurse cells, the border cell cluster, and the oocyte—each represented by their own phase-field variables \( \phi_i \) for \( i = 1, \dots, N \). {\allblack We
define the local cellular occupancy by
\begin{align*}
    \phi_{\max}(\mathbf{x},t)
    =
    \max_{1\leq i\leq N}\phi_i(\mathbf{x},t).
\end{align*}
The diffusion coefficient is then given by
\begin{align}\label{eq pfm diff}
    D(\boldsymbol{\phi})
    =
    D_0\,
    h\left(
    \frac{1}{
    1+\exp\left[(\phi_{\max}-\phi_*)/s\right]}
    \right),
\end{align}}
% We define the diffusion coefficient as:
% \begin{align}\label{eq pfm diff}
%     D(\phi) =D_0 h\left(\frac{1}{1 + \sum_{i=1}^N e^{(\phi_i - \phi_*)/s}}\right),
% \end{align}
where \( D_0 \) is the baseline diffusion constant, {\allblack $h(q) = q^2(3-2q)$} is the interpolation function, \( \phi_* \) {\allblack determines the phase-field
value at which the transition in diffusivity occurs,} 
%is the threshold parameter, 
and \( s \) {\allblack controls the steepness of this transition.}
% a steepness parameter controlling the transition in diffusivity. 
Eq.~\eqref{eq pfm diff} indicates that diffusion is reduced inside cellular regions (where \( \phi_i \approx 1 \)) and higher in extracellular regions (where all \( \phi_i \approx 0 \)).

The use of the interpolation function {\allblack \( h(q) = q^2(3 - 2q) \)} is critical for ensuring smooth and physically consistent transitions in the diffusion coefficient across cell boundaries. In the absence of interpolation function, \( D(\phi) \) would change sharply near phase-field interfaces, potentially introducing discontinuities in the spatial derivatives of the chemoattractant concentration. 
Such discontinuities could lead to numerical instabilities or spurious oscillations during simulation. 
The interpolation function $h$ regularizes sharp transitions so that the numerical scheme approximates smooth gradients instead of abrupt jumps.

\subsubsection{Coupling Tangential Interface Migration (TIM) Force to Spatial Chemoattractant Distribution}

Border cell migration in the \emph{Drosophila} egg chamber is shaped not only by chemoattractant gradients but also by physical contact with surrounding nurse cells~\cite{stonko2015mathematical, aranjuez2016dynamic, dai2020tissue, cai2014mechanical}. 
% {\allcpink need citations including Dai, maybe Arenjuez, Stonko, maybe others}. 
To represent this contact-mediated mechanism, we previously introduced the Tangential Interface Migration (TIM) force~\cite{akhavan2026phase}, which models traction generated along intercellular boundaries. This force reflects experimental observations of border cell clusters contacting and moving along adjacent nurse cells and generating protrusive forces at their interfaces~\cite{prasad2007cellular, aranjuez2016dynamic}.

In our earlier work, the TIM force was coupled to a prescribed or steady-state chemoattractant gradient \cite{george2025chemotaxis, akhavan2026phase}. In this study, we extend the model by allowing the chemoattractant field to evolve dynamically in space and time, thus capturing the effects of geometric constraints and spatial heterogeneity in signal distribution. 
This coupling provides a framework for investigating how the evolving chemoattractant distribution modulates mechanically driven interface migration through dynamic feedback.
We define the TIM force as:
\begin{align}\label{eq: tim force corrected}
\mathbf{F}_{\text{TIM}} = -\bar \mu_c \nabla \cdot \left(\rho(c) \, \phi_c \phi_j \, \operatorname{sgn}(\nabla c \cdot \nabla \phi_c^\perp) \left( \nabla \phi_c \right)^\perp \right),
\end{align}
where \( \phi_c \) and \( \phi_j \) represent the phase-field variables for the border cell cluster and neighboring nurse cells, respectively. 
% The function $\rho(c)$ represents receptor-mediated \cite{george2025chemotaxis, akhavan2026phase} sensing of the local chemoattractant concentration. 
{\allblack
The function $\rho(c)$ represents the effective receptor-mediated response to the local chemoattractant concentration and is defined as
\begin{align}
\rho(c)
=
\frac{s c^3}{(c^2+\Gamma)(c+\ell)},
\label{eq:rho}
\end{align}
where $s$ denotes the maximal activation level and $\Gamma$ and $\ell$ control the sensitivity of the response. This receptor-response function
is inherited from our previous chemoattractant--receptor model \cite{george2025chemotaxis,akhavan2026phase}.}
It captures the nonlinear dynamics of ligand–receptor interactions, including low sensitivity at minimal concentrations, heightened activation at intermediate levels, and saturation at high concentrations.
This sigmoidal response is consistent with experimentally observed PVF1–PVR signaling behavior in the egg chamber and was previously characterized in the modeling framework developed by George et al.~\cite{george2025chemotaxis} (see also~\cite{cai2016modeling, akhavan2026phase, dai2020tissue}). 
% {\allcpink need additional references for receptor dynamics here}  
The constant \( \bar \mu_c \) is the intensity parameter of TIM force, and $(\nabla\phi_c)^\perp$ denotes the tangential (perpendicular) {\allblack direction} of the cluster interface. 
The sign function $\operatorname{sgn}(\nabla c\cdot(\nabla\phi_c)^\perp)$ ensures that the force direction aligns with increasing chemoattractant along the interface.

{\allblack The use of the sign function is intentional in the TIM formulation. Unlike a conventional chemotactic force whose magnitude scales directly with the local chemoattractant gradient, TIM represents contact-mediated tangential traction that is activated according to the local receptor-response level $\rho(c)$. The projected gradient $\nabla c\cdot(\nabla\phi_c)^\perp$ therefore determines the preferred direction of tangential motion, while $\rho(c)$ modulates the magnitude of the chemically activated traction. Consequently, the present TIM formulation does not assume that the traction magnitude is proportional to the local gradient magnitude.

The use of $\operatorname{sgn}(\cdot)$ therefore represents an idealized directional switch. In regions where the projected chemoattractant gradient is extremely shallow, the present formulation does not account for a gradual loss of directional sensitivity. A smooth response such as $\tanh[\lambda\,\nabla c\cdot(\nabla\phi_c)^\perp]$ provides a natural extension for incorporating sensitivity to shallow-gradient magnitude.}

{\allblack
As in our previous phase-field formulation \cite{akhavan2026phase}, the migration force enters the right-hand side of the phase-field evolution equation for the border-cell cluster. In the TIM-driven simulations considered here, the conventional chemical force $F_{\mathrm{chem}}$ is replaced by $F_{\mathrm{TIM}}$, rather than being applied as an additional force.
Accordingly, the cluster phase field evolves as
\begin{align}
-\frac{1}{\mu}\frac{\partial \phi_c}{\partial t} &= \epsilon_c^2\nabla^2\phi_c + \phi_c(1-\phi_c) \left[ \phi_c-\frac{1}{2} +6F(\phi_c,\phi_0) \right] + F_{\mathrm{TIM}} .
\label{eq:phi_TIM}
\end{align}
% The remaining cellular phase fields evolve according to the phase-field dynamics described in Section~2 without this migration forcing term.
The remaining cellular phase fields evolve according to the phase-field
dynamics described under ``Phase Field Modeling of the Egg Chamber''
without this migration forcing term.

The TIM contribution is an active, non-variational forcing term and is not obtained from the free energy $E$. Consequently, after $F_{\mathrm{TIM}}$ is introduced, the cluster dynamics are no longer governed solely by a gradient flow, and monotonic decrease of $E$ is not guaranteed. Because $F_{\mathrm{TIM}}$ is written in divergence form, it does not introduce a net source of phase-field material when the corresponding
boundary flux vanishes. However, the cluster area $V_c=\int_\Omega h(\phi_c)\,d\mathbf{x}$ is not imposed as an exact
conservation law; deviations from its prescribed target value $\bar V_c$ are controlled by the existing volume-penalty term in the phase-field model.}

Eq.~\eqref{eq: tim force corrected} supports tangential migration along intercellular interfaces, with directionality influenced by both interface geometry and local chemoattractant gradients. By coupling the TIM force to a dynamically evolving chemoattractant field governed by a diffusion–reaction equation with geometry-dependent diffusivity, the model captures interactions between tissue structure, chemical signaling, and mechanical migration. This framework allows for the analysis of how spatial heterogeneity in chemoattractant availability, arising from extracellular architecture, influences the strength and direction of interface-driven migration.

\subsection{Numerical Implementation and Simulation Parameters}
We carry out simulations on a two-dimensional computational domain defined as \( \Omega = [0,5] \times [0,5] \), which is discretized uniformly into an \( M \times M \) grid. The spatial mesh size is set to \( h = \frac{5}{M+1} \), and temporal evolution is computed using a forward Euler scheme with fixed time step \( \Delta t \). 
{\allblack
The computational coordinates are expressed in model units. Associating the
modeled egg-chamber extent of 5 spatial units with approximately
$200\,\mu\mathrm{m}$ gives a characteristic length scale of
$40\,\mu\mathrm{m}$ per model unit. Accordingly, the spatial grid spacing
$\Delta x=\Delta y=0.05$ corresponds to approximately $2\,\mu\mathrm{m}$.}
% {\allblack Model parameters are selected within the phase-field framework to maintain stable cellular morphology, volume constraints, cell--cell interactions, and collective migration. All parameter values used in the simulations presented here are listed in Supplementary Table~1.}
{\allblack The phase-field formulation and the interpretation of its mechanical parameters follow our previous modeling framework~\cite{akhavan2026phase,nonomura2012study}. All parameter values used in the simulations presented here are listed in Supplementary Table~1.}
% Model parameters, including adhesion strength, volume penalization, and interfacial stiffness $(\alpha_m, \gamma(m,n), \beta(m,n))$, are selected to ensure biologically realistic morphology and collective behavior of the migrating cluster. Specific parameter values used in each simulation are listed in the corresponding figure captions \cite{nonomura2012study, akhavan2026phase}. 

The governing equations, including the phase-field dynamics and the diffusion–reaction equation for chemoattractant concentration, are solved using a finite difference method \cite{leveque2007finite} with Neumann (zero-flux) boundary conditions imposed on all edges of the domain \cite{boettinger2002phase, folch2005quantitative, mekus2018effects}. 
{\allblack For the present simulations, homogeneous Neumann boundary conditions are
used as a simplifying approximation to restrict chemoattractant transport
to the modeled egg-chamber cross-section. This condition does not imply that
the biological follicular epithelium is perfectly impermeable; rather, it
neglects protein exchange across the outer tissue boundary over the spatial and
temporal scales considered here.}
% These boundary conditions indicate the assumption that the egg chamber is surrounded by a larger tissue environment that does not permit chemoattractant or cell material to exit the computational domain. 
% All simulations are implemented in MATLAB R2022a, using a grid resolution of \( h = 0.05 \) and a time step of \( \Delta t = 0.05 \).
{\allblack
All simulations were implemented in MATLAB R2022a (The MathWorks, Inc., Natick, MA, USA) \cite{MATLAB2022a}, using a grid resolution of $h=0.05$ and a time step of $\Delta t=0.05$.
}

\section*{Data Availability}

No experimental datasets were used for model fitting or quantitative analysis in this study. The numerical results presented in this work were generated from the mathematical and computational models described in the manuscript and can be reproduced using the code provided in the Code Availability section. The original microscopy images supporting Fig.~10 are available from the corresponding author upon reasonable request.

\section*{Code Availability}
The computational code used to perform the simulations presented in this study is publicly available at \url{https://github.com/Naghmeh-Akhavan/PFM_Feedback_Concentration_Model}.

\section*{Acknowledgments}
% We acknowledge funding from NIH/NCCIH grant R01AT013188 to N.A., NSF DMS \#1953423 to B.E.P., and M.S.-G. Also, NSF IOS \#2303587 to M.S.-G. 
Not applicable.

\section*{Funding}
We acknowledge funding from NIH/NCCIH grant R01AT013188 to N.A.,
NSF DMS \#1953423 to B.E.P. and M.S.-G.,
NSF IOS \#2303587 to M.S.-G.,
and NSF DMS \#2619580 to N.A., A.K.A., M.S.-G., and B.E.P.

\section*{Author Contributions}
N.A. \& B.E.P. created the mathematical model to inform and guide the experimental framework. N.A. \& A.K.A. developed and carried out simulations and analyzed results from the model while A.S.G. carried out biological experiments and curated the resulting data.  All authors contributed to data interpretation, and reviewing and editing the text and visuals in the manuscript. The project was funded and overseen by M.S.-G. \& B.E.P.

\section*{Competing Interests}
%The authors declare no competing interests. 
The authors declare no competing financial or non-financial interests.

\bibliographystyle{naturemag}
\bibliography{References}

\newpage

\end{document}